\documentclass{amsart}
\usepackage{amssymb}
\usepackage{amscd}

\title[$p$-adic multiple zeta values I]
{$p$-adic multiple zeta values I \\
--- $p$-adic multiple polylogarithms 
and the $p$-adic KZ equation}
\author{Hidekazu Furusho}

\address{Research Institute for Mathematical Sciences\\ 
Kyoto University, Kyoto 606-8502 Japan}
\email{furusho@kurims.kyoto-u.ac.jp}

\newtheorem{thm}{Theorem}[section]
\newtheorem{lem}[thm]{Lemma}
\newtheorem{cor}[thm]{Corollary}
\newtheorem{prop}[thm]{Proposition}  

\theoremstyle{remark}
\newtheorem{ack}{Acknowledgments}        

\theoremstyle{definition}
\newtheorem{defn}[thm]{Definition}
\newtheorem{rem}[thm]{Remark}
\newtheorem{nota}[thm]{Notation}     
\newtheorem{notes}[thm]{Note}
\newtheorem{eg}[thm]{Examples}       
\newtheorem{pf}{Proof}                 
\newtheorem{q}[thm]{Question}   
\newtheorem{assump}[thm]{Assumption}

\numberwithin{equation}{section}

\begin{document}
\bibliographystyle{amsalpha+}
\maketitle

\begin{abstract}
Our main aim in this paper is to give a foundation of the theory of
$p$-adic multiple zeta values.
We introduce (one variable) $p$-adic multiple polylogarithms 
by Coleman's $p$-adic iterated integration theory.
We define $p$-adic multiple zeta values to be special values of
$p$-adic multiple polylogarithms.
We consider the (formal) $p$-adic KZ equation and introduce the
$p$-adic Drinfel'd associator
by using certain two fundamental solutions of the $p$-adic KZ equation.
We show that our $p$-adic multiple polylogarithms appear 
as coefficients of 
a certain fundamental solution of the $p$-adic KZ equation 
and our $p$-adic multiple zeta values appear 
as coefficients of the $p$-adic Drinfel'd associator.
We show various properties of $p$-adic multiple zeta values,
which are sometimes analogous to the complex case
and are sometimes peculiar to the $p$-adic case,
via the $p$-adic KZ equation.
\end{abstract}

\tableofcontents

\setcounter{section}{-1}
\section{Introduction}\label{winter}
The aim of the present paper and the upcoming papers \cite{F2} and \cite{F3}
is to enlighten crystalline aspects of
the fundamental group of the projective line minus three points and
add crystalline part to \cite{F1}.
In this paper, we will introduce the notions of 
(one-variable) $p$-adic multiple polylogarithms,
$p$-adic multiple zeta values, 
$p$-adic KZ equation and
$p$-adic Drinfel'd associator,
which will be our basic foundations of \cite{F2} and \cite{F3},
and show their various properties and their relationships.\par
Let $k_1,\cdots,k_m\in\bold N$.
The (usual) {\sf multiple zeta value} is the real number defined by the 
following series
\begin{equation}\label{Rokko}
\zeta(k_1,\cdots,k_m)=
\underset{n_i\in\bold  N}
{
\underset{0<n_1<\dotsm<n_m}{\sum}
}
\frac{1}{n_1^{k_1}\dotsm n_m^{k_m}}
\qquad  .
\end{equation}
Especially in the case when $m=1$, 
the multiple zeta value coincides with
the Riemann zeta value $\zeta(k)$.
We can check easily that this series converges in the topology of $\bold R$
if and only if $k_m> 1$,
however, this series never converges in the topology of $\bold Q_p$!
Thus it is not so easy and not so straightforward to give a definition of 
$p$-adic version of multiple zeta value.
To give a nice definition, we need another interpretation of 
multiple zeta values.\par
Suppose $z\in\bold C$.
The (one variable)
{\sf multiple polylogarithm} is a function defined by the following series
\begin{equation*}
Li_{k_1,\cdots,k_m}(z)=
\underset{n_i\in\bold  N}
{
\underset{0<n_1<\dotsm<n_m}{\sum}
}
\frac{z^{n_m}}{n_1^{k_1}\dotsm n_m^{k_m}}
\qquad  .
\end{equation*}
Especially in the case when $m=1$,
the multiple polylogarithm coincides with 
the classical polylogarithm $Li_k(z)$.
Easily we see that this series converges for $|z|<1$.
In \S\ref{Monday}, we will define the {\sf $p$-adic multiple polylogarithm} 
to be the function 
defined by the above series just replacing $z\in\bold C$ by $z\in\bold C_p$.
We remark that especially in the case when $m=1$,
the $p$-adic multiple polylogarithm is equal to 
the $p$-adic polylogarithm $\ell_k(z)$
which was studied by Coleman \cite{C}.
What is interesting is that
this $p$-adic multiple polylogarithm converges for $|z|_p<1$ similarly to
the above complex case.
Here $|\cdot|_p$ means the standard multiplicative valuation 
of $\bold C_p$.\par
An important relationship between the multiple polylogarithm and 
the multiple zeta value is the following formula:
\begin{equation}\label{Mikasa}
\zeta(k_1,\cdots,k_m)=
\underset{|z|<1}{
\underset{z\to 1}{\lim}}
Li_{k_1,\cdots,k_m}(z)
\qquad  .
\end{equation}
In this paper, we will define the $p$-adic multiple zeta value
\`a la formula \eqref{Mikasa} instead of \`a la \eqref{Rokko}.
But we  note that here happens a serious problem because
the open unit disk centered at $0$ on $\bold C_p$ and
the one centered at $1$ on $\bold C_p$ are disjoint!
Thus we cannot consider $\underset{z\to 1}{\lim}$ of
$p$-adic multiple polylogarithms
which are functions defined on $|z|_p<1$,
i.e. on the open unit disk centered at $0$.
To give a meaning of this limit, we will make an analytic continuation
of $p$-adic multiple polylogarithms 
by Coleman's $p$-adic iterated integration theory \cite{C}
and then define $p$-adic multiple zeta values
to be a limit value at $1$ of 
analytically continued $p$-adic multiple polylogarithms.\par

The organization of this paper is as follows.
\S\ref{spring} is devoted to a short review of well-known results on (usual)
multiple polylogarithms and multiple zeta values 
and definitions of the KZ equation and the Drinfel'd associator,
which will play a role of prototype in the $p$-adic case 
in the following two sections.\par
In \S\ref{summer}, we will introduce $p$-adic multiple zeta values and
show their many nice properties.
At first, we will review Coleman's 
$p$-adic iterated integration theory \cite{C}
in \S\ref{Sunday}
and then in \S\ref{Monday}
we will give an analytic continuation of 
$p$-adic multiple polylogarithms 
(which is just a multiple analogue of that of 
his $p$-adic polylogarithms $\ell_k(z)$ in \cite{C})
to the whole plane minus $1$, i.e. $\bold C_p-\{1\}$,
by his integration theory.
But we will see that there happens a terrible problem that 
the analytically continued $p$-adic multiple polylogarithm admits too many
(uncountably infinite) branches $Li_{k_1,\cdots,k_m}^a(z)$ 
$(k_1,\cdots,k_m\in\bold N, z\in\bold C_p-\{1\})$
which correspond to each branch parameter $a\in \bold C_p$,
coming from branch $log^a(z)$ of $p$-adic logarithms
(see \S\ref{Monday}).
However the following theorem in \S\ref{Tuesday} will remove our anxiety.\\
{\bf Theorem \ref{Gunma}.}
{\it
If $\underset{z\in\bold C_p-\{1\}}{\underset{z\to 1}{\lim}^\prime} Li_{k_1,\cdots,k_m}^a(z)$ converges on $\bold C_p$,
its limit does not depend on any choice of branch parameter $a\in \bold C_p$.
For $\lim^\prime$, see Notation \ref{Tochigi}.
}\\
By this theorem,
we can give a definition of the {\sf $p$-adic multiple zeta value}
$\zeta_p(k_1,\cdots,k_m)$ to be the above limit on $\bold C_p$
as follows.\\
{\bf Definition \ref{Saitama}.}
{\it
$\zeta_p(k_1,\cdots,k_m):=
\underset{z\in\bold C_p-\{1\}}{\underset{z\to 1}{\lim}^\prime} Li_{k_1,\cdots,k_m}^a(z)\in \bold C_p$
if it converges.
}\\
This definition of $p$-adic multiple zeta value is 
actually independent of any choice of branch parameter $a\in\bold C_p$
by Theorem \ref{Gunma}.
Especially, in the case when $m=1$,
we shall see in Example \ref{Tokyo} (due to Coleman \cite{C})
that the $p$-adic multiple zeta value is equal to
the $p$-adic $L$-value up to a certain constant multiple.
The following three theorems in \S\ref{Tuesday} 
are $p$-adic analogues of basic properties,
Lemma \ref{May}, 
Lemma \ref{August} and 
Proposition \ref{October},
in the complex case.\\
{\bf Theorem \ref{Chiba}.}
{\it
$\underset{z\in\bold C_p-\{1\}}{\underset{z\to 1}{\lim}^\prime} Li_{k_1,\cdots,k_m}^a(z)$
converges on $\bold C_p$ if $k_m>1$.
}\\
{\bf Theorem \ref{Ishikawa}.}
{\it
$\zeta_p(k_1,\cdots,k_m)\in \bold Q_p$.
}\\
{\bf Theorem \ref{Nagano}.}
{\it
The product of two $p$-adic multiple zeta values
can be written as a $\bold Q$-linear combination of 
$p$-adic multiple zeta values.
}\par
In \S\ref{autumn}, we will consider the $p$-adic KZ equation and introduce
the $p$-adic Drinfel'd associator which will play a role of main tools to prove
Theorem \ref{Chiba}, Theorem \ref{Niigata} and Theorem \ref{Nagano}.
In \S\ref{jazz}, we will introduce the {\sf $p$-adic KZ equation}
in Definition \ref{Palestine} and prove the following.\\
{\bf Theorem \ref{Egypt}.}
{\it
Let $a\in\bold C_p$. 
Then there exists a unique solution $G_0^a(A,B)(z)$
($z\in\bold P^1(\bold C_p)\backslash\{0,1,\infty\}$)
of the $p$-adic KZ equation which is a formal power series whose coefficients
are Coleman functions with respect to $a\in\bold C_p$ and
are locally analytic on
$\bold P^1(\bold C_p)\backslash\{0,1,\infty\}$
and satisfies a certain asymptotic behavior
$G_0^a(z)\approx z^A$ ($z\to 0$),
where $z^A:=1+\frac{log^a(z)}{1!}A+\frac{(log^a(z))^2}{2!}A^2+\cdots$.
}\\
Then we will  introduce a definition of {\sf $p$-adic Drinfel'd associator}
$\Phi^p_{KZ}(A,B)$ from two fundamental solutions,
$G_0^a(u)$ and $G_1^a(u)$, 
of the $p$-adic KZ equation in Definition \ref{Spain} 
and show its branch independency in Theorem \ref{Ireland}.
In \S\ref{waltz}, we will state precisely and prove the following.\\
{\bf Theorem \ref{Belgium}.}
{\it
Let $a\in\bold C_p$.
The fundamental solution $G_0^a(z)$ 
($z\in\bold P^1(\bold C_p)\backslash\{0,1,\infty\}$)
of the $p$-adic KZ equation can be expressed in terms of
(analytically continued) $p$-adic multiple polylogarithms 
$Li_{k_1,\cdots,k_m}^a(z)$ and the $p$-adic logarithm $log^a(z)$
explicitly.
}\\
In \S\ref{pop}, we will state precisely and prove the following.\\
{\bf Theorem \ref{Brazil}.}
{\it
The $p$-adic Drinfel'd associator $\Phi^p_{KZ}(A,B)$ 
can be expressed explicitly in terms of $p$-adic multiple zeta values.
}\\
In \S\ref{rock}, we will show functional equations of 
$p$-adic multiple polylogarithms at first.\\
{\bf Theorem \ref{Malaysia}.}
\[
Li_{k_1,\cdots,k_m}^a(1-z)=(-1)^m\underset{W=W'W''}{\sum_{W',W'':\text{words}}}
I_p(W'')\cdot J_p^a\bigl(\tau(W')\bigr)(z) \ 
\]
{\it
where $W=A^{k_m-1}BA^{k_{m-1}-1}B\cdots A^{k_1-1}B$.
Here each $I_p(W'')$ is a certain $\bold Q$-linear combination of 
$p$-adic multiple zeta values (see Theorem \ref{Brazil})
and each $J_p^a\bigl(\tau(W')\bigr)(z)$ is a certain combination of 
$p$-adic multiple polylogarithms (see Theorem \ref{Belgium}).
}\\
Next we will see that especially 
Coleman-Sinnott's functional equation of the
$p$-adic dilogarithm will be re-proved in Example \ref{Indonesia} 
by Theorem \ref{Malaysia}
and then will prove Theorem \ref{Niigata} and Theorem \ref{Nagano}.\par

In the upcoming paper \cite{F2}, 
we will relate the $p$-adic Drinfel'd associator with
the crystalline Frobenius action on the
rigid (unipotent) fundamental group of the projective line
minus three points 
and compare the $p$-adic Drinfel'd associator with other corresponding objects
in various realizations of motivic fundamental groups 
of the projective line minus three points. 
In \cite{F3}, certain algebraic relations 
among $p$-adic multiple zeta values are shown.
We will show there that the $p$-adic Drinfel'd associator determines a point
of the Grothendieck-Teichm\"{u}ller group
and discuss the elements of Ihara's stable derivation algebra 
arising from this point and 
finally we will add crystalline part to \cite{F1}.

\begin{ack}
The author wish to thank his advisor Professor Akio Tamagawa
for his encouragement.
The author is grateful to Hiroshi Fujiwara who helped the author to
type this manuscript.
The author is supported in part by 
JSPS Research Fellowships for Young Scientists.
\end{ack}

\section{Review of the complex case}\label{spring}
We shall review definitions of multiple polylogarithms and multiple zeta values
in \S\ref{day} and shall recall notions of the KZ equation and the 
Drinfel'd associator briefly in \S\ref{night},
which may help to understand its $p$-adic version developed 
in the following two sections.

\subsection{Multiple polylogarithms and multiple zeta values}\label{day}
We review briefly definitions and properties of 
multiple polylogarithms and multiple zeta values.
For more details, consult \cite{F0} and \cite{Gon} for example.\par
Let $k_1,\cdots,k_m\in\bold N$ and $z\in\bold C$.

\begin{defn}\label{January}
The (one variable) {\sf multiple polylogarithm} ({\sf MPL} for short) 
is defined to be the following series:
\begin{equation*}
Li_{k_1,\cdots,k_m}(z)=
\underset{n_i\in\bold  N}
{
\underset{0<n_1<\dotsm<n_m}{\sum}
}
\frac{z^{n_m}}{n_1^{k_1}\dotsm n_m^{k_m}}
\qquad  .
\end{equation*}
\end{defn}

\begin{rem}\label{New Year Day}
This MPL is the special case of the multiple polylogarithm
\[
Li_{k_1,\cdots,k_m}(z_1,\cdots,z_m)=
\underset{n_i\in\bold  N}
{
\underset{0<n_1<\dotsm<n_m}{\sum}
}
\frac{z_1^{n_1}\dotsm z_m^{n_m}}{n_1^{k_1}\dotsm n_m^{k_m}}
\]
introduced in \cite{Gon}
where $z_1=\cdots=z_{m-1}=1$ and $z_m=z$
\end{rem}

Easily we can check the following.

\begin{lem}\label{February}
The MPL $Li_{k_1,\cdots,k_m}(z)$ converges for $|z|<1$.
\end{lem}

\begin{lem}\label{March}
Suppose that $|z|<1$. Then
\begin{align*}
&\frac{d}{dz}Li_{k_1,\cdots,k_m}(z)=
\begin{cases}
\frac{1}{z}Li_{k_1,\cdots,k_m-1}(z) & k_m\neq 1, \\
\frac{1}{1-z}Li_{k_1,\cdots,k_{m-1}}(z)     & k_m=1, \\
\end{cases}  \\
&\frac{d}{dz}Li_1(z)=\frac{1}{1-z} \qquad .
\end{align*}
\end{lem}

By Lemma\ref{March}, for $|z|<1$, we get $Li_1(z)=-log(1-z)$
and the following,
\[
Li_{k_1,\cdots,k_m}(z)=
\begin{cases}
\int_0^z\frac{1}{t}Li_{k_1,\cdots,k_m-1}(t)dt & k_m\neq 1, \\
\int_0^z\frac{1}{1-t}Li_{k_1,\cdots,k_{m-1}}(t)dt & k_m= 1, \\
\end{cases}  
\]
from which  we get an expression of the MPL 
by iterated integral of $\frac{dt}{t}$ and $\frac{dt}{1-t}$.
Since  $\frac{dt}{t}$ and $\frac{dt}{1-t}$ admit poles at $t=0,1$ and $\infty$,
we cannot give an analytic continuation of the MPL to the whole complex plane 
due to monodromies around $0$, $1$ and $\infty$.
However we can say that

\begin{lem}\label{April}
The MPL $Li_{k_1,\cdots,k_m}(z)$ can be analytically continued to the universal
unramified covering $\widetilde{ \bold P^1(\bold C)\backslash\{0,1,\infty\} }$
of $\bold P^1(\bold C)\backslash\{0,1,\infty\}$.
\end{lem}

Since the simply-connected Riemann surface 
$\widetilde{ \bold P^1(\bold C)\backslash\{0,1,\infty\} }$
is an infinite covering of  $\bold P^1(\bold C)\backslash\{0,1,\infty\}$,
each MPL admits ({\em countably}) infinite branches.
The following are well-known ( for example, see \cite{Gon}).

\begin{lem}\label{May}
$\underset{|z|<1}{\underset{z\to 1}{\lim}}Li_{k_1,\cdots,k_m}(z)$
converges if $k_m>1$.
\end{lem}

\begin{lem}\label{June}
$\underset{|z|<1}{\underset{z\to 1}{\lim}}Li_{k_1,\cdots,k_m}(z)$
diverges if $k_m=1$.
\end{lem}

\begin{defn}\label{July}
For $k_1,\cdots,k_m\in \bold N$, $k_m>1$,
the {\sf multiple zeta value} ({\sf MZV} for short) is defined to be 
\[
\zeta(k_1,\cdots,k_m)=
\underset{|z|<1}{
\underset{z\to 1}{\lim}}
Li_{k_1,\cdots,k_m}(z)
\Biggl(=
\underset{0<n_1<\dotsm<n_m}{\sum}
\frac{1}{n_1^{k_1}\dotsm n_m^{k_m}}
\Biggr)
\qquad  .
\]
\end{defn}

Since MPL's are $\bold C$-valued functions, MZV's lie in $\bold C$.
However we can say more.

\begin{lem}\label{August}
For $k_1,\cdots,k_m\in \bold N$, $k_m>1$,
$\zeta(k_1,\cdots,k_m)\in\bold R$.
\end{lem}

\begin{nota}\label{September}
For each natural number $w$,
let $Z_w$ be the $\bold Q$-vector subspace of $\bold R$ generated by all MZV's 
$\zeta(k_1,\cdots,k_m)$ with $k_1+\cdots+k_m=w$, i.e.
$Z_w:=\langle\zeta(k_1,\cdots,k_m)\ | \  k_1+\cdots+k_m=w\rangle_{\bold Q}\subseteq\bold R$,
and put $Z_0=\bold Q$.
Denote $Z\centerdot$ to be the {\it formal direct sum} of 
$Z_w$ for all $w\geqslant 0$:
$Z\centerdot:=\underset{w\geqslant 0}{\oplus}Z_w$. 
\end{nota}

The following is one of the fundamental properties of MZV's.

\begin{prop}\label{October}
The graded $\bold Q$-vector space  $Z\centerdot$ forms 
a graded $\bold Q$-algebra, i.e.
$Z_a\cdot Z_b\subseteq Z_{a+b}$ for $a,b\geqslant 0$.
\end{prop}

\begin{pf}
At least we have two proofs \cite{Gon}.
The first one is given by the {\sf harmonic product formulae},
from which it follows, for example,
\[
\zeta(m)\cdot\zeta(n)=\zeta(m,n)+\zeta(n,m)+\zeta(m+n).
\]
The other one is given by the {\sf shuffle product formulae},
from which it follows, for example,
\[
\zeta(m)\cdot\zeta(n)=\sum_{i=0}^{m-1}\binom{n-1+i}{i}\zeta(m-i,n+i)+
\sum_{j=0}^{n-1}\binom{m-1+j}{j}\zeta(n-j,m+j).
\]
\qed
\end{pf}

\subsection{The KZ equation and the Drinfel'd associator}\label{night}
In this subsection, we will briefly review the definition of 
the KZ equation and the Drinfel'd associator.
For more detailed information on the KZ equation and the Drinfel'd associator,
see \cite{Dr}, \cite{F0} and \cite{Kas}.\par
Let $A^{\land}_{\bold C}=\bold C\langle\langle A,B\rangle\rangle$
be the non-commutative formal power series ring 
generated by two elements $A$ and $B$
with complex number coefficients.

\begin{defn}\label{November}
The (formal) {\sf Knizhnik-Zamolodchikov equation}
\footnote{
This is a special case of the KZ equation 
for $\bold P^1(\bold C)\backslash\{0,1,\infty\}$
in \cite{Kas}.
}
({\sf KZ equation} for short) is the differential equation
\begin{equation}
\tag*{(KZ)}
\frac{\partial G}{\partial u}(u)=
(\frac{A}{u}+\frac{B}{u-1})\cdot G(u) \ \ \ \ \ \ \ \ \ \ ,
\end{equation}
where $G(u)$ is an analytic function in complex variable $u$
with values in  $A^{\land}_{\bold C}$
where \lq analytic' means  
each of whose coefficient is analytic.
\end{defn}
The equation (KZ) has singularities only at $0,1$ and $\infty$.
Let $\bold C'$ be the complement of the union of the real 
half-lines $(-\infty,0]$ and $[1,+\infty)$ in the complex plane.
This is a simply-connected domain.
The equation (KZ) has a unique analytic solution on $\bold C'$ having a specified value 
at any given points on  $\bold C'$.
Moreover, for the singular points $0$ and $1$,
there exist unique solutions $G_0(u)$ and $G_1(u)$ of (KZ) such that
\[
G_0(u)\approx u^A  \  (u\to 0), \ \   \ \ \ \  G_1(u)\approx (1-u)^B  \  (u\to 1),  
\]
where $\approx$ means that 
$G_0(u)\cdot u^{-A}$ (resp. $G_1(u)\cdot(1-u)^{-B}$) 
has an analytic  continuation in a neighborhood of $0$ (resp. $1$) 
with value $1$ at $0$ (resp. $1$).
Here,
$u^A:=exp(A log u)
:=1+\frac{A log u}{1!}+\frac{(A log u)^2}{2!}+\frac{(A log u)^3}{3!}+\dotsb$
and
$log u := \int_1^u\frac{dt}{t}$ in $\bold C'$.
In the same way, $(1-u)^B$ is well-defined on $\bold C'$.
Since $G_0(u)$ and $G_1(u)$ are both non-zero unique solutions of (KZ) 
with the specified asymptotic behaviors, 
they must coincide with each other up to multiplication
from the right by an invertible element of $A^{\land}_{\bold C}$.

\begin{defn}\label{December}
The {\sf Drinfel'd associator}
\footnote{
To be precise, Drinfel'd defined $\varphi_{KZ}(A,B)$ instead of $\Phi_{KZ}(A,B)$ in \cite{Dr}, 
where $\varphi_{KZ}(A,B)=\Phi_{KZ}(\frac{1}{2\pi i}A,\frac{1}{2\pi i}B)$.}
is the  element $\Phi_{KZ}(A,B)$
of $A^{\land}_{\bold C}$ which is defined  by
\[
\ \ \ \ G_0(u)=G_1(u)\cdot\Phi_{KZ}(A,B) \  \ .
\]
\end{defn}

By considering on $(A^{\land}_{\bold C})^{ab}$, the abelianization of
$A^{\land}_{\bold C}$, we easily find that $\Phi_{KZ}(A,B)\equiv 1$ on
$(A^{\land}_{\bold C})^{ab}$.
We note that MZV's appear at each coefficient of the Drinfel'd associator
$\Phi_{KZ}(A,B)$.
For its explicit formulae, see \cite{F0}Proposition 3.2.3.

\section{$p$-adic multiple polylogarithms and $p$-adic multiple zeta values}
\label{summer}
In this section, we shall give the definition of 
$p$-adic multiple polylogarithms
(\S\ref{Monday}) and $p$-adic multiple zeta values (\S\ref{Tuesday})
and state main results in \S\ref{Tuesday},
which will be proved in \S \ref{autumn}.
The reader will find interesting analogies between \S\ref{day} and
\S\ref{Tuesday}.

\subsection{Review of Coleman's $p$-adic iterated integration theory}
\label{Sunday}
We will review the $p$-adic iterated integration theory 
by R. Coleman \cite{C},
following A. Besser's reformulation in \cite{Be1}.
This theory will be employed in the analytic continuation of 
$p$-adic multiple polylogarithms in \S\ref{Tuesday}.
For other nice expositions of his theory,
see \cite{Be2}\S 5, \cite{Br}\S 2.2.1 and \cite{CdS}\S 2.\par

\begin{assump}\label{Hakodate}
Suppose that $X/{\mathcal O}_{\bold C_p}$ is a smooth projective 
and surjective scheme over the ring ${\mathcal O}_{\bold C_p}$ 
of integers of $\bold C_p$,
of relative dimension $1$ 
with its generic fiber $X_{\bold C_p}$ and its special fiber
$X_{\overline{\bold F_p}}$.
Let $Y=X-D$ where $D$ is a closed subscheme of $X$
which is relatively etale over ${\mathcal O}_{\bold C_p}$.
\end{assump}

We denote 
$j:Y_{\overline{\bold F_p}}\hookrightarrow X_{\overline{\bold F_p}}$ 
to be the associated open embedding,
where $Y_{\overline{\bold F_p}}$ is the special fiber of $Y$ and
denote the finite set $X(\overline{\bold F_p})-Y(\overline{\bold F_p})$
by $\{e_1,\cdots,e_s\}$.
For $0\leqslant r < 1$,
$U_r$ stands for the rigid analytic space
\footnote{
As is explained in \cite{Be1} \S 2,
while the definition of $U_r$ depends on the choice of  `local lifts',
the definitions of $A^a_{log}(U_x)$ and $ \Omega^a_{log}(U_x)$
(see below) do not.
} 
obtained by 
{\it removing all closed discs of radius $r$ around $e_i$ from 
$X({\bold C_p})$} 
($1\leqslant i \leqslant s$) (see \cite{Be1}).
For a subset $S\subset X(\overline{\bold F_p})$,
we denote its tubular neighborhood (see \cite{Ber})
in $X({\bold C_p})$ by $]S[$.
For any rigid analytic space $W$, we mean by $A(W)$ the ring of 
global sections of the sheaf ${\mathcal O}^{\mathrm{ rig}}_W$ of
rigid analytic functions on $W$.\par

Fix $a\in\bold C_p$.
It determines a branch of $p$-adic logarithm 
$log^a:\bold C_p^\times\to\bold C_p$ (\cite{Be1}Definition 2.6) 
which is characterized by $log^a(p)=a$.
We call this $a\in\bold C_p$ the {\sf branch parameter}
of $p$-adic logarithm. 
Define
\[
A^a_{\mathrm{loc}}:=\prod_{x\in X(\overline{\bold F_p})}A^a_{log}(U_x),\qquad
\Omega^a_{\mathrm{loc}}:=\prod_{x\in X(\overline{\bold F_p})}\Omega^a_{log}(U_x)
\]
where
\begin{align*}
&A^a_{log}(U_x):=
\begin{cases}
A(]x[) &  x\in Y(\overline{\bold F_p}),\\
\underset{r\to 1}{\mathrm{ind\text{-}lim }}
A\Bigl(]x[ \ \cap U_r\Bigr)\Bigl[log^a(z_x)\Bigr] &
x\in\{e_1,\cdots,e_s\},
\end{cases}\\
&\Omega^a_{log}(U_x):=A^a_{log}(U_x)dz_x \ \ .
\end{align*}
Here $z_x$ means a local parameter 
\[
z_x: \ ]x[ \ \cap Y({\bold C_p})\overset{\sim}{\to}\{z\in\bold C_p \ |  \ 0<|z|_p<1\}.
\]
We note that $log^a(z_x)$ is a locally analytic function defined on
$]x[ \ \cap Y({\bold C_p})$
whose derivation is $\frac{1}{z_x}$ and it is transcendental over 
$\underset{r\to 1}{\mathrm{ind\text{-}lim }}A\Bigl(]x[ \ \cap U_r\Bigr)$
and $\underset{r\to 1}{\mathrm{ind\text{-}lim }}A\Bigl(]x[ \ \cap U_r\Bigr)\cong\Bigl\{f(z)=\sum\limits_{n=-\infty}^{n=\infty}a_nz^n  \ (a_n\in\bold C_p)
\text{ converging for }r<|z|_p<1 \text{ for some } 0\leqslant r<1\Bigr\}$
(see \cite{Be1}).
We remark that these definitions of $A^a_{log}(U_x)$ and $\Omega^a_{log}(U_x)$
are independent of any choice of local parameters $z_x$.
By taking a component-wise derivative, we obtain a $\bold C_p$-linear map
$d:A^a_{\mathrm{loc}}\to\Omega^a_{\mathrm{loc}}$.
Regard $A^\dag:=\Gamma(]X_{\overline{\bold F_p}}[, j^\dag\mathcal{O}_{]X_{\overline{\bold F_p}}[})$ and
$\Omega^\dag:=\Gamma(]X_{\overline{\bold F_p}}[, j^\dag\Omega^1_{]X_{\overline{\bold F_p}}[})$
to be a subspace of $A^a_{\mathrm{loc}}$ and $\Omega^a_{\mathrm{loc}}$ 
respectively
(for $j^\dag$, see \cite{Ber}).\par
In \cite{C}, Coleman constructed an 
$A^\dag$-subalgebra $A^a_{\mathrm{Col}}$ of $A^a_{\mathrm{loc}}$,
which we call {\sf the ring of Coleman functions} 
attached to a branch parameter $a\in\bold C_p$,
and a $\bold C_p$-linear map
$\int_{(a)}:A^a_{\mathrm{Col}}\underset{A^\dag}{\otimes}\Omega^\dag\to A^a_{\mathrm{Col}}\Bigm/ \bold C_p\cdot 1$
satisfying $d\bigm|_{A^a_{\mathrm{Col}}}\circ\int_{(a)}=id_{A^a_{\mathrm{Col}}\otimes\Omega^\dag}$,
which we call {\sf the $p$-adic (Coleman) integration}
attached to a branch parameter $a\in\bold C_p$.
We often drop  the subscript ${}_{(a)}$. \par

Actually Coleman's $p$-adic integration theory is essentially independent 
of any choice of branches, which may not be well-known, and
we will try to explain this fact below:\\
Suppose that $a,b\in\bold C_p$. 
Consider the isomorphisms
$\iota_{a,b}:A^a_{\mathrm{loc}}\overset{\sim}{\to}A^b_{\mathrm{loc}}$ and
$\tau_{a,b}:\Omega^a_{\mathrm{loc}}\overset{\sim}{\to}\Omega^b_{\mathrm{loc}}$
obtained by replacing each $log^a(z_{e_i})$ by $log^b(z_{e_i})$
for $1\leqslant i\leqslant s$.

\begin{lem}\label{Sapporo}
These maps, $\iota_{a,b}$ and $\tau_{a,b}$,
are independent of any choice of 
a local parameter $z_{e_i}$.
\end{lem}

\begin{pf}
Suppose that $z'_{e_i}$ is another local parameter.
Then we check easily that 
$\iota_{a,b}(z'_{e_i})=z'_{e_i}$,
$\iota_{a,b}(log^a(z'_{e_i}))=log^b(z'_{e_i})$ and
$\tau_{a,b}(dz'_{e_i})=dz'_{e_i}$
because $log^a(\frac{z'_{e_i}}{z_{e_i}})$ is analytic at $e_i$,
from which it follows the lemma.
\qed
\end{pf}

The following branch independency principle,
which was not stated explicitly in \cite{C},
should be one of the important properties of 
Coleman's $p$-adic integration theory,
but this principle just follows directly from his construction of 
$A^a_{\mathrm{Col}}$.

\begin{prop}[Branch Independency Principle]\label{Hokkaido}
Suppose that $a,b\in\bold C_p$. Then
$\iota_{a,b}(A^a_{\mathrm{Col}})=A^b_{\mathrm{Col}}$,
$\tau_{a,b}(A^a_{\mathrm{Col}}\otimes\Omega^\dag)=A^b_{\mathrm{Col}}\otimes\Omega^\dag$ and
$\iota_{a,b}\circ\int_{(a)}=\int_{(b)}\circ\tau_{a,b} \mod \bold C_p\cdot 1$.
Namely the following diagram is commutative.
\[
\begin{CD}
A^a_{\mathrm{Col}}\underset{A^\dag}{\otimes}\Omega^\dag
@>{\tau_{a,b}}>>
A^b_{\mathrm{Col}}\underset{A^\dag}{\otimes}\Omega^\dag \\
@V{\int_{(a)}}VV
@VV{\int_{(b)}}V \\
A^a_{\mathrm{Col}}\Bigm/ \bold C_p\cdot 1
@>{\iota_{a,b}}>>
A^b_{\mathrm{Col}}\Bigm/ \bold C_p\cdot 1 \\
\end{CD}
\]
\end{prop}

\begin{pf}
It follows directly because both $(A^b_{\mathrm{Col}},\int_{(b)})$
and 
$(\iota_{a,b}(A^b_{\mathrm{Col}}),
\iota_{a,b}\circ\int_{(a)}\circ\tau_{a,b}^{-1})$
satisfies the same axioms (A)--(F)
of logarithmic $F$-crystals on \cite{C}
(see also \cite{CdS}),
which is  a characterization of 
$(A^b_{\mathrm{Col}},\int_{(b)})$.
\qed
\end{pf}

Other important properties of Coleman's functions are the uniqueness principle 
and the functorial property below.

\begin{prop}[Uniqueness Principle ; \cite{C}Ch IV]\label{Aomori}
Let $a\in\bold C_p$.
Let $f\in A^a_{\mathrm{Col}}$ be a Coleman function which is defined on
an admissible open subset $U$ of $X(\bold C_p)$.
Suppose that $f|_U\equiv 0$.
Then $f\equiv 0$ on $X(\bold C_p)$.
\end{prop}

Especially this proposition yields the fact that a locally constant 
Coleman function is globally constant.

\begin{prop}
[Functorial Property ; \cite{C}Theorem 5.11 and \cite{Be2}Definition 4.7]
\label{Iwate}
Let $a\in\bold C_p$.
Let $(X',Y')$ be another pair satisfying 
Assumption \ref{Hakodate}.
Suppose that $f:X'\to X$ is a morphism defined over $\mathcal{O}_{\bold C_p}$
such that $f(Y')\subset Y$.
Then the pull-back morphism 
$f^\sharp:A^a_{\mathrm{loc}}\to A^{\prime a}_{\mathrm{loc}}$
induces the morphism
$f^*:A^a_{\mathrm{Col}}\to A^{\prime a}_{\mathrm{Col}}$
of rings of Coleman functions,
where $A^{\prime a}_{\mathrm{loc}}$ (resp. $A^{\prime a}_{\mathrm{Col}}$)
stands for  $A^a_{\mathrm{loc}}$ (resp. $A^a_{\mathrm{Col}}$) for $(X',Y')$.
\end{prop}

Precisely speaking, Coleman showed this theorem in more general situation in
\cite{C} Theorem 5.1.

\begin{nota}\label{Miyagi}
Let $a\in\bold C_p$ and 
$\omega\in A^a_{\mathrm{Col}}\underset{A^\dag}{\otimes}\Omega^\dag$.
Then by Coleman's integration theory, there exists
(uniquely modulo constant) a Coleman function $F_\omega$ such that
$\int\omega\equiv F_\omega$ (modulo constant).
For $x,y\in]Y(\overline{\bold F_p})[$, 
we define $\int_x^y\omega$ to be $F_\omega(y)-F_\omega(x)$.
It is clear that $\int_x^y\omega$ does not depend on any choice
of $F_\omega$ (although it may depend on $a\in\bold C_p$).
If $F_\omega(x)$ and $F_\omega(y)$ make sense for some $x,y\in X(\bold C_p)$,
we also denote $F_\omega(y)-F_\omega(x)$ by $\int_x^y\omega$.
When we let $y$ vary, we regard $\int_x^y\omega$ as the Coleman function
which is characterized by $dF_\omega=\omega$ and $F_\omega(x)=0$.
\end{nota}

\subsection{Analytic continuation of $p$-adic multiple polylogarithms}
\label{Monday}
We will define $p$-adic multiple polylogarithms to be Coleman functions 
which admits an expansion around $0$ similar to the complex case.\par
Let $k_1,\cdots,k_m\in\bold N$ and $z\in\bold C_p$.
Consider the following series 
\begin{equation}\label{Fujisan}
Li_{k_1,\cdots,k_m}(z)=
\underset{n_i\in\bold  N}
{
\underset{0<n_1<\dotsm<n_m}{\sum}
}
\frac{z^{n_m}}{n_1^{k_1}\dotsm n_m^{k_m}}
\qquad  .
\end{equation}

\begin{lem}\label{Akita}
This series $Li_{k_1,\cdots,k_m}(z)$ converges on the open unit disk
$D(0:1)=\{ z\in\bold C_p \ | \ |z|_p<1\}$ around $0$ with radius $1$. 
\end{lem}

\begin{pf}
Easy.
\qed
\end{pf}

\begin{lem}\label{Yamagata}
Let $z\in\bold C_p$ such that $|z|_p<1$. Then
\begin{align*}
&\frac{d}{dz}Li_{k_1,\cdots,k_m}(z)=
\begin{cases}
\frac{1}{z}Li_{k_1,\cdots,k_m-1}(z) & k_m\neq 1, \\
\frac{1}{1-z}Li_{k_1,\cdots,k_{m-1}}(z)     & k_m=1, \\
\end{cases}  \\
&\frac{d}{dz}Li_1(z)=\frac{1}{1-z} \qquad .
\end{align*}
\end{lem}

\begin{pf}
It follows from a direct calculation.
\qed
\end{pf}

Lemma \ref{Akita} and Lemma \ref{Yamagata} are $p$-adic analogue
of Lemma \ref{February} and Lemma \ref{March} respectively.\par

From now on, we fix a branch parameter $a\in\bold C_p$ and employ Coleman's
$p$-adic integration theory 
attached to this branch parameter $a\in\bold C_p$
for $X=\bold P^1_{{\mathcal O}_{\bold C_p}}$ and
$Y=Spec \ {\mathcal O}_{\bold C_p}[t,\frac{1}{t},\frac{1}{1-t}]$.

\begin{defn}\label{Fukushima}
We define recursively the (one variable) {\sf $p$-adic multiple polylogarithm}
({\sf $p$-adic MPL} for short)
$Li^a_{k_1,\cdots,k_m}(z)\in A^a_{\mathrm{Col}}$
attached to $a\in\bold C_p$
which is the Coleman function characterized below:
\begin{align*}
& Li^a_{k_1,\cdots,k_m}(z):=
\begin{cases}
\int_0^z\frac{1}{t}Li^a_{k_1,\cdots,k_m-1}(t)dt & k_m\neq 1, \\
\int_0^z\frac{1}{1-t}Li^a_{k_1,\cdots,k_{m-1}}(t)dt & k_m= 1, \\
\end{cases}  \\
& Li^a_1(z)=-log^a(1-z):=\int_0^z\frac{dt}{1-t} \quad .
\end{align*}
\end{defn}  

\begin{rem}\label{Iwashiro}
\begin{enumerate}
\renewcommand{\labelenumi}{(\arabic{enumi})}
\item
It is easy to see that 
$Li^a_{k_1,\cdots,k_m}(z)=\underset{0<n_1<\dotsm<n_m}{\sum}
\frac{z^{n_m}}{n_1^{k_1}\dotsm n_m^{k_m}}$ if $|z|_p<1$.
\item
Because we get $\frac{dt}{t}$, $\frac{dt}{1-t}\in A^a_{\mathrm{Col}}$
and we know that $Li^a_{k_1,\cdots,k_m}(z)$ is analytic on $|z|_p<1$
and takes value $0$ at $z=0$,
it is easy to see that each $p$-adic MPL is well defined in
$A^a_{\mathrm{Col}}$.
\item
Our construction of $p$-adic MPL  is just a multiple analogue
of Coleman's construction \cite{C} of $p$-adic polylogarithm
$\ell_k(z)$.
His $p$-adic polylogarithm $\ell_k(z)$
can be written as $Li_k^a(z)$ in our notation.
\end{enumerate}
\end{rem}

The following is a $p$-adic version of Lemma \ref{April}.

\begin{prop}\label{Ibaraki}
The $p$-adic MPL $Li^a_{k_1,\cdots,k_m}(t)$ is locally analytic on
$\bold P^1(\bold C_p)\backslash\{1,\infty\}$.
More precisely,
$Li^a_{k_1,\cdots,k_m}(t)\bigm|_{]0[}\in A(]0[)$,
$Li^a_{k_1,\cdots,k_m}(t)\bigm|_{]1[}\in A(]1[)\bigl[log^a(t-1)\bigr]$ and
$Li^a_{k_1,\cdots,k_m}(t)\bigm|_{]\infty[}\in A(]\infty[)\bigl[log^a(\frac{1}{t})\bigr]$.
\end{prop}

\begin{pf}
By construction, we can prove the claim inductively.
\qed
\end{pf}

This proposition means that the series \eqref{Fujisan} can be analytically continued to
$\bold P^1(\bold C_p)\backslash\{1,\infty\}$,
although the complex MPL cannot be analytically continued 
to $\bold P^1(\bold C)\backslash\{0,1,\infty\}$
but only to its universal unramified covering 
$\widetilde{ \bold P^1(\bold C)\backslash\{0,1,\infty\} }$
instead by Lemma \ref{April}.
We note that the $p$-adic MPL admits {\it uncountably} infinite branches
which correspond to branches of $p$-adic logarithms $log^a(z)$
although complex MPL admits {\it countably} infinite branches.
We call $Li^a_{k_1,\cdots,k_m}(z)$ 
{\sf the branch of $p$-adic MPL} corresponding to $a\in\bold C_p$.

\subsection{$p$-adic multiple zeta values and main results}\label{Tuesday}
We will state main results of this paper,
whose proof will be given in \S\ref{pop} and \S\ref{rock}
and will introduce $p$-adic multiple zeta values,
whose definitions themselves are highly non-trivial.

\begin{nota}\label{Tochigi}
Let $\alpha\in\bold C_p$ and let $f(z)$ be a function
defined on $\bold C_p$.
We denote 
$\underset{z\to\alpha}{\lim^\prime} f(z)$ 
to be
$\underset{n\to\infty}{\lim} f(z_n)$
if this limit converges to the same value for any sequence
$\{z_n\}_{n=1}^{\infty}$ which satisfies
$z_n\to\alpha$ in $\bold C_p$ and 
$e(\bold Q_p(z_1,z_2,\cdots)/\bold Q_p)<\infty$
(which means that the field generated by $z_1,z_2,\cdots$ over $\bold Q_p$ is
a finitely ramified (possibly infinite) extension field over $\bold Q_p$).
If the latter limit converges (resp. does not converge) to the same value,
we call $\underset{z\to\alpha}{\lim^\prime} f(z)$ 
converges (resp. diverges).
\end{nota}

\begin{thm}\label{Gunma}
Fix $k_1,\cdots,k_m\in\bold N$ and a prime $p$.
Then the statement whether 
$\underset{z\in\bold C_p-\{1\}}{\underset{z\to 1}{\lim}^\prime} Li_{k_1,\cdots,k_m}^a(z)$ converges or diverges
\footnote{
This makes sense because $p$-adic MPL's are locally analytic on
$\bold P^1(\bold C_p)\backslash\{1,\infty\}$
by Proposition \ref{Ibaraki}.}
on $\bold C_p$ is independent of any choice of branch parameter
$a\in\bold C_p$.
Moreover if it converges on $\bold C_p$,
this limit value is independent of any choice of 
branch parameter $a\in\bold C_p$.
\end{thm}

\begin{pf}
Since $Li_{k_1,\cdots,k_m}^a(z)$ is locally analytic on 
$\bold P^1(\bold C_p)\backslash\{1,\infty\}$
by Proposition \ref{Ibaraki},
its $A^a_{log}(]1[)$-component can be written as follows.
\begin{align*}
Li_{k_1,\cdots,k_m}^a(z)=&f_0(z-1)+f_1(z-1)log^a(z-1)+f_2(z-1)(log^a(z-1))^2\\
& +\cdots\cdots +f_m(z-1)(log^a(z-1))^m,
\end{align*}
where $f_i(z)\in A(D(0:1))$ for $i=0,\cdots,m$.
By Proposition \ref{Hokkaido}, we see that these $f_i(z)$'s are independent 
of any choice of branch parameter $a\in\bold C_p$.
Saying 
$\underset{z\to 1}{\lim}^\prime Li_{k_1,\cdots,k_m}^a(z)$ converges is equivalent to saying $f_i(0)=0$
for all $i=1,\cdots,m$ by Lemma \ref{Russia} and Lemma \ref{Maebashi},
which is a statement independent of any choice of branch parameter 
$a\in\bold C_p$.
Thus we get the first half of this theorem.
If $\underset{z\to 1}{\lim}^\prime Li_{k_1,\cdots,k_m}^a(z)$ converges,
then $\underset{z\to 1}{\lim}^\prime Li_{k_1,\cdots,k_m}^a(z)=f_0(0)$.
Since $f_0(z)$ was independent of any choice of branches, 
the second half of this theorem follows.
\qed
\end{pf}

\begin{lem}\label{Russia}
For $n\geqslant 0$, $\underset{\epsilon\in\bold C_p}{\underset{\epsilon\to 0}{\lim}^\prime} \epsilon (log^a\epsilon )^n=0$.
\end{lem}

\begin{pf}
If $n=0$, it is clear.\\
If $n=1$, suppose that $\epsilon_n\in L$ ($n\geqslant 1$) and
$\epsilon_n\to 0$ as $n\to\infty$,
where $L$ is a finitely ramified (possibly infinite) extension of 
$\bold Q_p$ with ramification index $e_L(<\infty)$ and a uniformizer $\pi_L$.
Take $c_L\in\bold N$ such that $p^{c_L}>e_L$.
Decompose $\epsilon_n=u_n\cdot\pi_L^{r_n}$
where $u_n\in{\mathcal O}_L^\times$ and $r_n\in\bold Z$.
Take $s_n\in\bold N$ such that $(s_n,p)=1$ and 
$u_n^{s_n}\equiv 1 \mod \pi_L\mathcal O_L$.
Put $\alpha_n:=u_n^{s_n}-1\in\pi_L\mathcal O_L$.
Then $(u_n^{s_n})^{p^{c_L}}=(1+\alpha_n)^{p^{c_L}}\equiv 1+\alpha_n^{p^{c_L}}\equiv 1 \mod p\mathcal O_{\bold C_p}$.
Therefore $log^a u_n=\frac{1}{s_n\cdot p^{c_L}}log^a(u_n^{s_n\cdot p^{c_L}})\in\frac{1}{p^{c_L}}\mathcal O_{\bold C_p}$.
So we get 
\[
\lim_{n\to\infty}\epsilon_n log^a\epsilon_n=
\lim_{n\to\infty}\{\epsilon_n log^a u_n+\epsilon_n r_n log^a\pi_L\}
=\lim_{n\to\infty}\epsilon_n log^a u_n=0.
\]
In a similar way, we can prove the case for $n>1$.
\qed
\end{pf}

\begin{lem}\label{Maebashi}
Let $a\in\bold C_p$, $l\geqslant 0$ and 
$g(z)=\sum\limits_{k=0}^{l}a_k\bigl(log^a(z)\bigr)^k$
($a_k\in\bold C_p$).
Then $\underset{z\to 1}{\lim}^\prime g(z)$ converges
if and only if $a_k=0$ for $1\leqslant k \leqslant l$.
\end{lem}

\begin{pf}
Take $z_n=\alpha^n$  such that $|\alpha|_p<1$
and $log^a(\alpha)\neq 0$.
Then we get the claim by an easy calculation.
\qed
\end{pf}

\begin{rem}\label{maebashi}
\begin{enumerate}
\renewcommand{\labelenumi}{(\arabic{enumi})}
\item
As for another proof of the last statement of Theorem \ref{Gunma}
for $k_m>1$, see Remark \ref{Cuba}.
\item
It is striking that 
$\underset{\epsilon\to 0}{\lim}^\prime Li_{k_1,\cdots,k_m}^a(1-\epsilon)$
does not depend on any choice of branch parameter $a\in\bold C_p$
although $Li_{k_1,\cdots,k_m}^a(1-\epsilon)$ takes whole vales on $\bold C_p$
if we fix $\epsilon$ ($0<|\epsilon|_p<1$) and let $a$ vary on $\bold C_p$.
\end{enumerate}
\end{rem}

\begin{defn}\label{Saitama}
For any index $(k_1,\cdots,k_m)$ whose 
$\underset{z\to 1}{\lim}^\prime Li_{k_1,\cdots,k_m}^a(z)$ converges,
we define the corresponding {\sf $p$-adic multiple zeta vale}
({\sf $p$-adic MZV} for short) $\zeta_p(k_1,\cdots,k_m)$
to be its limit in $\bold C_p$, i.e.
\[
\zeta_p(k_1,\cdots,k_m):=
\underset{z\in\bold C_p-\{1\}}{\underset{z\to 1}{\lim}^\prime} Li_{k_1,\cdots,k_m}^a(z)\in\bold C_p\qquad
\text{       if it converges.}
\]
If this limit diverges for $(k_1,\cdots,k_m)$,
we do not give a definition of the corresponding $p$-adic MZV
$\zeta_p(k_1,\cdots,k_m)$.
\end{defn}

We note that this definition of $p$-adic MZV is independent of any choice of
branches by
Theorem \ref{Gunma}.

\begin{thm}\label{Chiba}
If $k_m>1$, $\underset{z\in\bold C_p-\{1\}}{\underset{z\to 1}{\lim}^\prime} Li_{k_1,\cdots,k_m}^a(z)$ always converges on $\bold C_p$.
\end{thm}

Thus we get a definition of $p$-adic MZV $\zeta_p(k_1,\cdots,k_m)$ 
for $k_m>1$.
This theorem is a $p$-adic analogue of Lemma \ref{May}
and it will be proved in \S\ref{pop}.

\begin{eg}\label{Tokyo}
Coleman made the following calculation in \cite{C}
(stated in Ch I (4) and proved in Ch VII):
\begin{equation}\label{Ikoma}
\underset{z\to 1}{\lim}^\prime Li_n^a(z)=
\frac{p^n}{p^n-1}L_p(n,\omega^{1-n})\quad \text{for  } n>1.
\end{equation}
Here $L_p$ is the Kubota-Leopoldt $p$-adic $L$-function and
$\omega$ is the Teichm\"{u}ller character.
In particular this formula \eqref{Ikoma} shows that
this limit value is actually independent of any choice of branch parameter 
$a\in\bold C_p$
although this fact is a special case of Theorem \ref{Gunma}.
We remark that, in the case of $p$-adic polylogarithms, 
this branch independency also follows from the so-called 
distribution relation (\cite{C} Proposition 6.1).
By \eqref{Ikoma}, we get
\begin{equation}\label{Hida}
\zeta_p(n)=\frac{p^n}{p^n-1}L_p(n,\omega^{1-n})\quad \text{for  } n>1.
\end{equation}
\begin{enumerate}
\renewcommand{\labelenumi}{(\alph{enumi})}
\item 
When $n$ is even (i.e. $n=2k$ for some $k\geqslant 1$),
by \eqref{Hida} we get the equality
\[
\zeta_p(2k)=0.
\]
We will see that this equality proved arithmetically here
will be also deduced from geometric identities, 
2-cycle relation and 3-cycle relation, among $p$-adic MZV's,
proved in \cite{F3}.
\item
On the other hand, when $n$ is odd (i.e. $n=2k+1$ for some $k\geqslant 0$),
it does not look so easy to show $\zeta_p(2k+1)\neq 0$:\\
Suppose that $p$ is an odd prime. 
Then by \cite{KNQ} Theorem 3.1,
we see that saying $L_p(2k+1,\omega^{-2k})\neq 0$
is equivalent to saying 
\begin{equation}\tag{L$_{2k+1}$}\label{Haiffa}
H^2_{\mathrm{et}}(\bold Z,\bold Q_p/\bold Z_p(-k))=0 ,
\end{equation}
which is one of standard conjectures 
\footnote{
In \cite{KNQ}, (L$_{2k+1}$) was denoted by (C$_k$).
}
in Iwasawa theory and a higher version 
of Leopoldt conjecture (cf. \cite{KNQ} Remark 3.2.(ii)).
\end{enumerate}
\end{eg}

\begin{rem}\label{Ganei-Yeelim}
\begin{enumerate}
\renewcommand{\labelenumi}{(\roman{enumi})}
\item
We know that $\zeta_p(2k+1)\neq 0$ in the case
where $p$ is regular or $(p-1)|2k$ by \cite{Sou}\S 3.3.
\item
Suppose that $p$ is an odd prime.
Let $G$ denote the standard Iwasawa module for 
$\bold Q(\mu_{p^\infty})\slash\bold Q$ \ 
($\mu_{p^\infty}$: the group of roots of unity whose order is a power of $p$),
i.e. $G=\pi_1\left(\bold Z[\mu_{p^\infty},\frac{1}{p}]\right)^{p,\text{ab}}$.
Let $I$ denote the inertia subgroup
of the unique prime $\frak p$ in $\bold Q(\mu_{p^\infty})$
which is above $p$.
Soul\'{e} \cite{Sou} constructed a specific non-zero element
$\chi_{p,m}^{\text{Soul\'{e}}}\in Hom\left(G,\bold Z_p(m)\right)$
for $m\geqslant 1$: odd.
Let $\mathcal G$ denote the standard Iwasawa module for 
$\bold Q_p(\mu_{p^\infty})\slash\bold Q_p$,
i.e. $\mathcal G=\pi_1\left(\bold Q_p(\mu_{p^\infty})\right)^{p,\text{ab}}$.
Let $\mathcal I$ denote the inertia subgroup of $\mathcal G$.
For $m\geqslant 1$, the Coates-Wiles homomorphism gives
a specific non-zero element 
$\chi_{p,m}^{\text{CW}}\in
Hom\left(\mathcal I,\bold Z_p(m)\right)$.
Coleman---using his reciprocity law---showed the following formula
for all odd $m>1$ (see \cite{KNQ}):
\begin{equation}\label{Tel-Aviv}
\chi_{p,m}^{\text{Soul\'{e}}}\circ r =(p^{m-1}-1)\cdot
L_p(m, \omega^{1-m})\cdot \chi_{p,m}^{\text{CW}}.
\end{equation}
Here $r$ is the natural (surjective) map $r:\mathcal I \to I$.
By \eqref{Tel-Aviv}, we get the following statement 
in algebraic number theory 
which is equivalent to saying $\zeta_p(2k+1)\neq 0$
(or, equivalently, \eqref{Haiffa} ):
\begin{equation}\tag{P$_{2k+1}$}
\text{
the prime ideal $\frak p$ ramify in the kernel field 
of $\chi_{p,m}^{\text{Soul\'{e}}}$.
}
\end{equation}
The author guesses more generally that problems on $p$-adic MZV's related to 
$p$-adic transcendental number theory
(such as the problem of proving the $p$-adic version 
$\zeta_p(3)\not\in\bold Q$ of Ap\'{e}ry's result)
could be translated into problems in algebraic number theory.
\end{enumerate}
\end{rem}

As for a $p$-adic analogue of Lemma \ref{June}, at present,
we have nothing to say except the following.

\begin{notes}\label{Kanagawa}
The limit $\underset{z\in\bold C_p-\{1\}}{\underset{z\to 1}{\lim}^\prime} Li_{k_1,\cdots,k_m}^a(z)$
sometimes converges and sometimes diverges on $\bold C_p$ for $k_m=1$.
\end{notes}

For example, see Example \ref{Toyama}.(a) and (b) below.

\begin{thm}\label{Niigata}
Suppose that $\underset{z\in\bold C_p-\{1\}}{\underset{z\to 1}{\lim}^\prime} Li_{k_1,\cdots,k_m}^a(z)$ converges on $\bold C_p$ for $k_m=1$.
Then it converges to a $p$-adic version of the regularized MZV, i.e.
\[
\zeta_p(k_1,\cdots,k_{m-1},1)=(-1)^mI_p(W)
\text{   where   }
W=BA^{k_{m-1}-1}B\cdots A^{k_1-1}B.
\]
\end{thm}

See \S\ref{waltz} for $I_p(W)$ and Remark \ref{Peru}(2) for the
regularized $p$-adic MZV.
This theorem will be proved in \S\ref{rock}.
Therefore $p$-adic MZV $\zeta_p(k_1,\cdots,k_m)$ for $k_m=1$
can be written as a $\bold Q$-linear combination of $p$-adic MZV's
corresponding to the same weight indexes with $k_m>1$.

\begin{eg}\label{Toyama}
\begin{enumerate}
\renewcommand{\labelenumi}{(\alph{enumi})}
\item 
$\underset{z\to 1}{\lim}^\prime Li_{2,1}^a(z)$
converges to $-2\zeta_p(1,2)$, i.e. $\zeta_p(2,1)=-2\zeta_p(1,2)$.
This follows from the functional equation in Example \ref{Indonesia}.(a),
$\zeta_p(2)=0$ by Example \ref{Tokyo}.(a) 
and Lemma \ref{Russia}.
\item
$\underset{z\to 1}{\lim}^\prime Li_{3,1}^a(z)$
diverges if and only if $\zeta_p(3)\neq 0$
(equivalently if and only if $p$ satisfies 
the 3rd Leopoldt conjecture (L$_3$) above).\par
Suppose that 3rd Leopoldt conjecture (L$_3$) fails at a prime $p$.
Then we get $\zeta_p(3,1)=-2\zeta_p(1,3)-\zeta_p(2,2)$
for this prime $p$.
This follows from the functional equation in Example \ref{Indonesia}.(b)
combined with Lemma \ref{Russia}.
\item
We will show many identities between $p$-adic MZV's in \cite{F3},
from which we will deduce, for example,
$\zeta_p(3)=\zeta_p(1,2)$ and $\zeta_p(1,3)=\zeta_p(2,2)=\zeta_p(1,1,2)=0$.
\end{enumerate}
\end{eg}

\begin{rem}\label{Viet-Nam}
The author guesses that to know whether 
$\underset{z\in\bold C_p-\{1\}}{\underset{z\to 1}{\lim}^\prime} Li_{k_1,\cdots,k_m}^a(z)$ converges or diverges 
might be to tell something deep in number theory,
such as Example \ref{Toyama}.(b).
\end{rem}

Those $p$-adic MZV's {\it were} defined to be elements of $\bold C_p$,
but actually we can say more.

\begin{thm}\label{Ishikawa}
All $p$-adic MZV's are $p$-adic numbers, i.e. 
$\zeta_p(k_1,\cdots,k_m)\in\bold Q_p$.
\end{thm}

\begin{pf}
Suppose that $\underset{z\in\bold C_p-\{1\}}{\underset{z\to 1}{\lim}^\prime} Li_{k_1,\cdots,k_m}^a(z)$ converges.
Recall that $p$-adic MPL  $Li^a_{k_1,\cdots,k_m}(z)$ $(a\in\bold C_p)$
is an iterated integral of $\frac{dt}{t}$ and $\frac{dt}{1-t}$
which is a rational $1$-form defined over $\bold Q_p$
and notice that  $Li^a_{k_1,\cdots,k_m}(z)\in\bold Q_p$ for all
$z\in p\bold Z_p$.
Then from the Galois equivariancy stated in \cite{BdJ} Remark 2.3,
it follows that $Li^a_{k_1,\cdots,k_m}(z)$ is 
$Gal(\overline{\bold Q_p}/\bold Q_p)$-invariant for
$z\in\bold P^1(\bold Q_p)\backslash\{1,\infty\}$
if we take $a\in\bold Q_p$.
Therefore in this case, we get $Li_{k_1,\cdots,k_m}^a(z)\in\bold Q_p$
for $z\in\bold P^1(\bold Q_p)\backslash\{1,\infty\}$.
Thus we get
$\underset{z\in\bold Q_p-\{1\}}{\underset{z\to 1}{\lim}} Li_{k_1,\cdots,k_m}^a(z)\in\bold Q_p$, which yields the theorem
(Recall that this limit is independent of any choice of branch parameter
$a\in\bold C_p$ by Theorem \ref{Gunma}).
\qed
\end{pf}

It may be better to say that this theorem is a $p$-adic version of
Lemma \ref{August}.
The author poses the following question, 
which he wants to study in the future.

\begin{q}\label{Fukui}
Are all $p$-adic MZV's $p$-adic integers? 
Namely $\zeta_p(k_1,\cdots,k_m)\in\bold Z_p$ for all primes $p$?
\end{q}

\begin{defn}\label{Yamanashi}
For each natural number $w$, let $Z_w^{(p)}$ be the finite dimensional
$\bold Q$-linear subspace of $\bold Q_p$ generated by all
$p$-adic MZV's of indices with weight $w$,
and put $Z_0^{(p)}=\bold Q$.
Define $Z^{(p)}_\centerdot$ to be the {\it formal direct sum}
of $Z_w^{(p)}$ for all $w\geqslant 0$:
$Z^{(p)}_\centerdot:=\underset{w\geqslant 0}{\oplus}Z^{(p)}_w$.
\end{defn}

By Theorem \ref{Niigata}, we see that 
$Z_w^{(p)}=\langle \zeta_p(k_1,\cdots,k_m) \ | \ k_1+\cdots+k_m=w,k_m>1,m\in\bold N\rangle_{\bold Q}\subset \bold Q_p$

\begin{thm}\label{Nagano}
The graded $\bold Q$-vector space has a structure of $\bold Q$-algebra,
i.e. $Z_a^{(p)}\cdot Z_b^{(p)}\subseteq Z_{a+b}^{(p)}$ for $a,b\geqslant 0$.
\end{thm}

This is a $p$-adic analogue of Proposition \ref{October},
whose proof will be given in \S \ref{rock}.
Unfortunately we do not have such a simple proof as Proposition \ref{October}.
Our proof is based on showing the {\sf shuffle product formulae}
(Corollary \ref{Singapore}) 
coming from the shuffle-like multiplication of iterated integrals,
from which it follows, for example,
\[
\zeta_p(m)\cdot\zeta_p(n)=\sum_{i=0}^{m-1}\binom{n-1+i}{i}\zeta_p(m-i,n+i)+
\sum_{j=0}^{n-1}\binom{m-1+j}{j}\zeta_p(n-j,m+j).
\]
In \cite{BF}, we shall discuss the {\sf harmonic product formulae}
\cite{H} coming from the shuffle-like multiplication of series in 
Definition \ref{January},
from which it should follow, for example,
\[
\zeta_p(m)\cdot\zeta_p(n)=\zeta_p(m,n)+\zeta_p(n,m)+\zeta_p(m+n).
\]
We note that the validity of the harmonic product formulae for
$p$-adic MZV's is non-trivial
because we do not have a series expansion of $p$-adic MZV
such as \S\ref{winter} \eqref{Rokko}.\par

\begin{rem}\label{Gifu}
Here is another direction of further possible developments of
our theory of $p$-adic MZV's.
Since the $p$-adic $L$-function is related to the Bernoulli numbers,
the author expects that {\sf the multiple Bernoulli number}
({\sf MBN} for short)
$B^{(k_1,\cdots,k_m)}_n\in\bold Q$ ($k_1,\cdots,k_m\in\bold Z$, $n\in\bold N$)
given by the following generating series
\[
\frac{Li_{k_1\cdots,k_m}(1-e^{-x})}{(1-e^{-x})^m}=:\sum_{n=0}^{\infty}
B^{(k_1,\cdots,k_m)}_n
\frac{x^n}{n!} \qquad ,
\]
where $Li_{k_1\cdots,k_m}(z)$ and $e^z$ means the following formal power series
\[
\underset{n_i\in\bold  N}{\underset{0<n_1<\dotsm<n_m}{\sum}}
\frac{z^{n_m}}{n_1^{k_1}\dotsm n_m^{k_m}}\in\bold Q[[z]]\quad
\text{ and }\quad
\sum_{n=0}^{\infty}\frac{z^n}{n!}\in\bold Q[[z]]
\]
respectively, would help to describe a $p$-adic behavior of 
$p$-adic MZV's and constructions of $p$-adic multiple zeta (or $L$-)functions.
We remark that this definition of MBN is just a multiple version of 
that of the poly-Bernoulli number in \cite{AK} and \cite{Kan}
(especially $B^1_n$ is the usual Bernoulli number).
We stress that the definition of MBN is independent of any prime $p$.
\end{rem}

\section{The $p$-adic KZ equation}\label{autumn}
In this section, we will introduce and consider the $p$-adic KZ equation.
We will give the definition of the $p$-adic Drinfel'd associator
$\Phi^p_{KZ}(A,B)$ in \S\ref{jazz}.
In \S\ref{waltz} (resp. \S\ref{pop}), 
we will give an explicit formula of a certain fundamental solution
$G^p_0(z)$ of the $p$-adic KZ equation 
(resp. an explicit formula of $\Phi^p_{KZ}(A,B)$).
In \S\ref{rock}, we will show the functional equation among $p$-adic MPL's
and will give proofs of Theorem \ref{Niigata} and Theorem \ref{Nagano}.

\subsection{The $p$-adic Drinfel'd associator}\label{jazz}
\begin{nota}\label{Yemen}
Let $A^{\land}_{\bold C_p}=\bold C_p\langle\langle A,B\rangle\rangle$
be the non-commutative formal power series ring with 
$\bold C_p$-coefficients generated by two elements $A$ and $B$.
\end{nota}

\begin{defn}\label{Palestine}
The (formal) {\sf $p$-adic Knizhnik-Zamolodchikov equation} 
({\sf $p$-adic KZ equation} for short) is the differential equation
\begin{equation}
\tag*{(KZ$^p$)}
\frac{\partial G}{\partial u}(u)=
(\frac{A}{u}+\frac{B}{u-1})\cdot G(u) \ \ \ \ \ \ \ \ \ \ ,
\end{equation}
where $G(u)$ is an analytic function in variable 
$u\in\bold P^1(\bold C_p)\backslash\{0,1,\infty\}$
with values in  $A^{\land}_{\bold C_p}$
where \lq analytic' means  
each of whose coefficient is locally $p$-adic analytic.
\end{defn}

Unfortunately, 
because $\bold P^1(\bold C_p)\backslash\{0,1,\infty\}$
is topologically totally disconnected,
the equation (KZ$^p$) does not have a unique solution on
$\bold P^1(\bold C_p)\backslash\{0,1,\infty\}$
even locally as in the complex analytic function case.
But fortunately we get the following nice property on Coleman functions.

\begin{thm}\label{Egypt}
Fix $a\in\bold C_p$. 
Then there exists a unique (invertible) solution 
$G_0^a(u)\in A^a_{\mathrm{Col}}\widehat\otimes A^{\land}_{\bold C_p}$
of (KZ$^p$) which is defined and locally analytic on
$\bold P^1(\bold C_p)\backslash\{0,1,\infty\}$ 
and satisfies $G_0^a(u)\approx u^A$ ($u\to 0$).
\end{thm}

Here $u^A:=1+\frac{log^a(u)}{1!}A+\frac{(log^a(u))^2}{2!}A^2+\cdots$
and $G_0^a(u)\approx u^A$ ($u\to 0$) means that the 
$A^a_{log}(]0[)\widehat\otimes A^{\land}_{\bold C_p}$-component of
$P^a(u):=G_0^a(u)\cdot u^{-A}=G_0^a(u)\cdot\Bigl\{1-\frac{log^a(u)}{1!}A+\frac{(log^a(u))^2}{2!}A^2-\frac{(log^a(u))^3}{3!}A^3+\cdots\Bigr\}$
lies in 
$1+A(]0[)\widehat\otimes A^{\land}_{\bold C_p}\cdot A
+A(]0[)\widehat\otimes A^{\land}_{\bold C_p}\cdot B$
and takes value $1\in A^{\land}_{\bold C_p}$ at $u=0$.

\begin{lem}\label{Turkey}
Fix $a\in\bold C_p$.
Let $G(u)$ and $H(u)$ be solutions of (KZ$^p$) in 
$A^a_{\mathrm{Col}}\widehat\otimes A^{\land}_{\bold C_p}$.
Suppose that $H(u)$ is invertible.
Then $H(u)^{-1}G(u)$ is a constant function,
i.e. an element of $A^{\land}_{\bold C_p}$.
\end{lem}

\begin{pf}
\begin{align*}
\frac{d}{du}& H(u)^{-1}G(u)=
-H(u)^{-1}\cdot\frac{d}{du}H(u)\cdot H(u)^{-1}G(u)+ H(u)^{-1}\frac{d}{du}G(u)\\
&=-H(u)^{-1}(\frac{A}{u}+\frac{B}{u-1})H(u)\cdot H(u)^{-1}G(u)+
H(u)^{-1}(\frac{A}{u}+\frac{B}{u-1})G(u)\\
&=-H(u)^{-1}(\frac{A}{u}+\frac{B}{u-1})G(u)+
H(u)^{-1}(\frac{A}{u}+\frac{B}{u-1})G(u)=0 \ .
\end{align*}
Therefore 
$H(u)^{-1}G(u)\in A^a_{\mathrm{Col}}\widehat\otimes A^{\land}_{\bold C_p}$ 
is constant by Proposition \ref{Aomori}.
\qed
\end{pf}

{\bf Proof of Theorem \ref{Egypt}.}\par
{\it Uniqueness}:
Suppose that $H_0^a(u)$ is another solution satisfying above properties.
Then $H_0^a(u)$ is invertible
(i.e. $H_0^a(u)\in (A^a_{\mathrm{Col}}\widehat\otimes A^{\land}_{\bold C_p})^\times$)
because it follows easily that
its constant term is $1\in A^a_{\mathrm{Col}}$.
By Lemma \ref{Turkey}, there exists a unique 'constant' series 
$c\in A^{\land}_{\bold C_p}$ such that 
$G_0^a(u)=H_0^a(u)\cdot c$.
By the assumption, we get $u^A\cdot c\cdot u^{-A}\to 1$ as $u\to 0$,
from which we can deduce $c=1$.
\par
{\it Existence}:
By substituting $G(u)=P(u)\cdot u^{A}$ into (KZ$^p$), we get 

\begin{equation}\label{Wakakusa}
\frac{dP}{du}(u)=\Bigl[\frac{A}{u},P(u)\Bigr ]+\frac{B}{u-1}P(u) \ \ .
\end{equation}

By expanding $P(u)=1+\sum\limits_{W:\text{words}} P_W(u)W $,
we obtain the following differential equation from \eqref{Wakakusa}:

\begin{align*}
&\frac{dP_W}{du}(u)=\frac{1}{u}P_{W'A}(u)-\frac{1}{u}P_{AW'}(u) 
&\text{if } W &=AW'A  \ (W'\in A^{\land}_{\bold C}),\\
&\frac{dP_W}{du}(u)=\frac{1}{u}P_{W'B}(u)
&\text{if } W &=AW'B  \ (W'\in A^{\land}_{\bold C}),\\
&\frac{dP_W}{du}(u)=\frac{-1}{u}P_{BW'}(u)+\frac{1}{u-1}P_{W'A}(u) 
&\text{if } W &=BW'A  \ (W'\in A^{\land}_{\bold C}),\\
&\frac{dP_W}{du}(u)=\frac{1}{u-1}P_{W'B}(u)
&\text{if } W &=BW'B  \ (W'\in A^{\land}_{\bold C}),\\
&\frac{dP_W}{du}(u)=0 
&\text{if } W &=A, \\
&\frac{dP_W}{du}(u)=\frac{1}{u-1}
&\text{if } W &=B. \\
\end{align*}
Since $\frac{du}{u}$ and $\frac{du}{1-u}$ lie in 
$A^a_{\mathrm{Col}}\otimes\Omega^\dag$
and are defined and locally analytic on 
$\bold P^1(\bold C_p)\backslash\{0,1,\infty\}$,
we can construct inductively a unique solution 
$P^a_0(u)=1+\sum\limits_{W:\text{words}} P^a_{0,W}(u)W $ of \eqref{Wakakusa}
such that each $P^a_{0,W}(u)$ satisfy 
the above differential equation,
lies in $A^a_{\mathrm{Col}}$,
is defined and locally analytic on 
$\bold P^1(\bold C_p)\backslash\{0,1,\infty\}$
and $P^a_W(0)=0$.
By putting $G_0^a(u)=P^a_0(u)\cdot u^A$,
we get a required solution in Theorem \ref{Egypt}.
\qed

\begin{prop}\label{Ethiopia}
Let $a,b\in\bold C_p$.
Then $\iota_{a,b}(G^a_0)=G^b_0$.
\end{prop}

\begin{pf}
This follows from the unique characterization of $G^b_0$ 
in Theorem \ref{Egypt}.
\qed
\end{pf}

\begin{prop}\label{Israel}
Fix $a\in\bold C_p$, $z_0\in\bold P^1(\bold C_p)\backslash\{0,1,\infty\}$
and $g_0\in A^{\land}_{\bold C_p}$.
Then there exists a unique solution 
$H^a(u)\in A^a_{\mathrm{Col}}\widehat\otimes A^{\land}_{\bold C_p}$
of (KZ$^p$) which satisfies $H^a(z_0)=g_0$.
Here $A^a_{\mathrm{Col}}\widehat\otimes A^{\land}_{\bold C_p}$
means the non-commutative two variable formal power series ring
with $A^a_{\mathrm{Col}}$-coefficients, i.e.
$A^a_{\mathrm{Col}}\widehat\otimes A^{\land}_{\bold C_p}=
A^a_{\mathrm{Col}}\langle\langle A,B\rangle\rangle$.
\end{prop}

\begin{pf}
This immediately follows from Lemma \ref{Turkey}
by taking $H^a(u)=G_0^a(u)\cdot G_0^a(z_0)^{-1}\cdot g_0$.
\qed
\end{pf}

\begin{prop}\label{Jordan}
Fix $a\in\bold C_p$.
Then there exists a unique solution 
$G_1^a(u)\in A^a_{\mathrm{Col}}\widehat\otimes A^{\land}_{\bold C_p}$
of (KZ$^p$) which is locally analytic on
$\bold P^1(\bold C_p)\backslash\{0,1,\infty\}$ and
satisfies $G_1^a(u)\approx (1-u)^B$ ($u\to 1$).
\end{prop}

Here the meanings of notations $(1-u)^B$ and 
$G_1^a(u)\approx (1-u)^B$ ($u\to 1$) are similar to those of 
$u^A$ and $G_0^a(u)\approx u^A$ ($u\to 0$) in Theorem \ref{Egypt}.

\begin{pf}
By a similar argument to Theorem \ref{Egypt},
we get the claim.
\qed
\end{pf}

\begin{prop}\label{England}
$G^a_1(u)=G^a_0(B,A)(1-u)$.
\end{prop}

Here for any $g\in A^{\land}_{\bold C_p}$, $g(B,A)$ stands for the image of
$g$ by the automorphism $A^{\land}_{\bold C_p}$
induced from $A\mapsto B$, $B\mapsto A$ and, for $f(u)\in A^a_{\mathrm{loc}}$,
$f(1-u)$ means its image by the algebra homomorphism
$\tau^\sharp: A^a_{\mathrm{loc}}\to A^a_{\mathrm{loc}}$
induced from the automorphism 
$\tau:t\mapsto 1-t$ of $\bold P^1(\bold C_p)\backslash\{0,1,\infty\}$.

\begin{pf}
By the functorial property of Coleman functions (Proposition \ref{Iwate}),
$\tau^\sharp(A^a_{\mathrm{Col}})\subseteq A^a_{\mathrm{Col}}$.
Therefore the proposition follows immediately from the uniqueness of
$G^a_1(u)$ because $G^a_0(B,A)(1-u)$ lies in 
$A^a_{\mathrm{Col}}\widehat\otimes A^{\land}_{\bold C_p}$,
satisfies (KZ$^p$) and admits the same asymptotic behavior to that of 
$G_1(u)$ at $u=1$.
\qed
\end{pf}

\begin{rem}\label{Britain}
As in the same way as the above proof of the uniqueness on Theorem \ref{Egypt},
we see that $G^a_0(u)$ and $G^a_1(u)$ are both invertible and 
they must coincide with each other up to a multiplication from the right 
by an invertible element 
$\Phi^{(a),p}_{KZ}(A,B)\in A^{\land}_{\bold C_p}=\bold C_p\langle\langle A,B\rangle\rangle$
(which is independent of $u\in\bold P^1(\bold C_p)\backslash\{0,1,\infty\}$).
Namely 
\[
G^a_0(u)=G^a_1(u)\cdot\Phi^{(a),p}_{KZ}(A,B).
\]
\end{rem}

\begin{thm}\label{Ireland}
Actually $\Phi^{(a),p}_{KZ}(A,B)$ is independent of any choice of branch
parameter $a\in\bold C_p$.
\end{thm}

\begin{pf}
Put 
$z_0\in\Bigl]\bold P^1_{\overline{\bold F_p}}\backslash\{0,1,\infty\}\Bigr[$.
Since $\Bigl]\bold P^1_{\overline{\bold F_p}}\backslash\{0,1,\infty\}\Bigr[$
is a branch independent region,
the special value of $G^a_0(u)$ at $u=z_0$
actually does not depend on any choice of branch parameter $a\in\bold C_p$.
Similarly we see that the value of $G^a_1(u)$ at $u=z_0$
does not depend on any choice of branch parameter.
Therefore $\Phi^{(a),p}_{KZ}(A,B)=G^a_1(z_0)^{-1}\cdot G^a_0(z_0)$
is actually independent of any choice of branch parameter $a\in\bold C_p$.
\qed
\end{pf}

\begin{rem}\label{France}
We have another proof of Theorem \ref{Ireland}.
Put $a,b\in\bold C_p$.
Then we get
$\iota_{a,b}(G^a_0(u))=G^b_0(u)$ and  $\iota_{a,b}(G^a_1(u))=G^b_1(u)$
by Proposition \ref{Ethiopia}.
Therefore we get $\iota_{a,b}(\Phi^{(a),p}_{KZ})=\Phi^{(b),p}_{KZ}$,
which implies $\Phi^{(a),p}_{KZ}=\Phi^{(b),p}_{KZ}$
because $\Phi^{(a),p}_{KZ}$ and $\Phi^{(b),p}_{KZ}\in A^{\land}_{\bold C_p}$.
\end{rem}

\begin{defn}\label{Spain}
The {\sf $p$-adic Drinfel'd associator} $\Phi^p_{KZ}(A,B)$ is 
the element of $\bold C_p\langle\langle A,B\rangle\rangle^\times$,
which is defined by $G_0^a(u)=G_1^a(u)\cdot\Phi^p_{KZ}(A,B)$.
\end{defn}

\begin{rem}\label{Portuguese}
\begin{enumerate}
\renewcommand{\labelenumi}{(\arabic{enumi})}
\item
This definition of 
the $p$-adic Drinfel'd associator  $\Phi^p_{KZ}(A,B)$ is independent 
of $u\in\bold P^1(\bold C_p)\backslash\{0,1,\infty\}$
by Remark \ref{Britain} and
any choice of branch parameter $a\in\bold C_p$
by Theorem \ref{Ireland}.
\item
In \S \ref{pop}, we will see that each coefficient of $\Phi^p_{KZ}(A,B)$
can be expressed in terms of $p$-adic MZV's,
from which  we know that actually $\Phi^p_{KZ}(A,B)$ belongs to
$\bold Q_p\langle\langle A,B\rangle\rangle$ by Theorem \ref{Ishikawa}.
\item
We shall prove many identities, such as $2$-, $3$- and $5$-cycle relation
of the $p$-adic Drinfel'd associator $\Phi^p_{KZ}(A,B)$
in \cite{F3}.
\end{enumerate}
\end{rem}

\subsection{Explicit formulae of the fundamental solution of the $p$-adic KZ equation}\label{waltz}
In this subsection, we will give a calculation to express each coefficient of
the fundamental solution $G^a_0(z)$ of the $p$-adic KZ equation (KZ$^p$).

\begin{nota}\label{Germany}
\begin{enumerate}
\renewcommand{\labelenumi}{(\arabic{enumi})}
\item
Let $\bold A\centerdot=\underset{w\geqslant 0}{\oplus}\bold A_w=\bold Q\langle A,B\rangle (\subset A^{\land}_{\bold C_p})$
be the non-commutative graded polynomial ring 
over $\bold Q$
with two variables $A$ and $B$ with
$deg A= deg B =1$.
Here $\bold A_w$ is the homogeneous degree $w$ part of $\bold A\centerdot$.
We call an element of $\bold A\centerdot$ which is a monomial 
with coefficient $1$
by a {\sf word}.
But exceptionally we shall not call $1$ a word.
For each word $W$, the {\sf weight} and {\sf depth} of $W$ are as follows.\par
$wt (W):=$ `the sum of exponents of $A$ and $B$ in $W$' \par
$dp (W):=$ `the sum of exponent  \ of\ \ \ \ \ \ \ \ \ \ $B$ in $W$'
\item
Put $M'=\bold A\centerdot\cdot B=\{ F\cdot B| F\in\bold A\centerdot\}$
which is the $\bold Q$-linear subspace of 
$\bold A\centerdot$.
Note that there is a natural surjection from
$\bold A\centerdot$  to
$\bold A\centerdot\Bigm/\bold A\centerdot A$.
By identifying the latter space with 
$\bold Q\cdot 1+M' (=\bold Q\cdot 1+\bold A\centerdot\cdot B)$
we obtain the $\bold Q$-linear map
$f':\bold A\centerdot\twoheadrightarrow\bold A\centerdot\Bigm/\bold A\centerdot A\overset{\sim}{\to}\bold Q\cdot 1+M'\hookrightarrow\bold A\centerdot$.
Abusively we denote by $f'$ the $\bold C_p$-linear map
$A^{\land}_{\bold C_p}\to A^{\land}_{\bold C_p}$
induced by
$f':\bold A\centerdot\to\bold A\centerdot$.
\item
For each word $W=B^{q_0}A^{p_1}B^{q_1}A^{p_2}B^{q_2}\dotsm A^{p_k}B^{q_k}$
($k\geqslant 0$, $q_0\geqslant 0$, $p_i,q_i\geqslant 1$ for $i\geqslant 1$)
in $M'$, we define
\[
Li^a_W(z):=Li^a_{
\underbrace{1,\ldots 1}_{q_k-1},
p_k+1,
\underbrace{1,\ldots 1}_{q_{k-1}-1},
p_{k-1}+1,\ldots\ldots ,1,p_1+1
\underbrace{1,\ldots 1}_{q_0}
}(z)\in A^a_{\mathrm{Col}} \ .
\]
By extending linearly, we get the $\bold Q$-linear map
$Li^a(z):M'\to A^a_{\mathrm{Col}}$
which sends each word $W$ in $M'$ to $Li^a_W(z)$.
\end{enumerate}
\end{nota}

\begin{thm}[Explicit Formulae]\label{Belgium} 
Put $a\in\bold C_p$.
Let $G^a_0(z)$ be the fundamental solution of the $p$-adic KZ equation 
(KZ$^p$) in Theorem \ref{Egypt}.
Expand $G^a_0(z)=1+\sum\limits_{W:\text{words}} J_p^a(W)(z) \ W $.
Then each coefficient $ J_p^a(W)(z)$ can be expressed as follows.
\begin{enumerate}
\renewcommand{\labelenumi}{(\alph{enumi})}
\item  When $W$ is in $M'$, \ \ \ $J_p^a(W)(z)=(-1)^{dp (W)}Li^a_W(z)$.
\item When $W$ is written as $VA^r (r\geqslant 0,V\in M')$,
\[
J_p^a(W)(z)=\underset{0\leqslant s,0\leqslant t}{\sum_{s+t=r}}
(-1)^{dp (W)+s}Li^a_{f'(V\circ A^s)}(z)
\frac{\{log^a(z)\}^t}{t!}\quad .
\]
\item When $W$ is written as 
$A^r (r\geqslant 0)$, \ \ \ $J_p^a(W)(z)=\frac{\{log^a(z)\}^r}{r!}$.
\end{enumerate}
\end{thm}

For the definition of the shuffle product \lq $\circ$',
see \cite{F0} Definition 3.2.2.
The proof of this theorem will be given in the end of this subsection.\par

\begin{lem}\label{Sweden}
$f'(G^a_0(z))=1+\sum\limits_{W\in M':\text{words}} J_p^a(W)(z) \ W $.
\end{lem}

\begin{pf}
Apply $f'$ term by term.
\qed
\end{pf}

\begin{lem}\label{Netherland}
$f'(G^a_0(z))=1+\sum\limits_{W\in M':\text{words}} (-1)^{dp(W)}Li_W^a(z) \ W $.
\end{lem}

\begin{pf}
By the $p$-adic KZ equation, we get the following:
\begin{align*}
&\frac{d}{du}J^a_p(W)(u)=\frac{1}{u}J^a_p(W')(u)
&\text{if } W &=AW' \ (W'\in M'),\\
&\frac{d}{du}J^a_p(W)(u)=\frac{1}{u-1}J^a_p(W')(u)
&\text{if } W &=BW' \ (W'\in M'),\\
&\frac{d}{du}J^a_p(W)(u)=\frac{1}{u-1}
&\text{if } W &=B\in M'.\\
\end{align*}
By Lemma \ref{Yamagata}, we see that the family
$\left\{(-1)^{dp(W)}Li_W^a(z)\in A^a_{\mathrm{Col}}\right\}_{W\in M':\text{words}}$
satisfies the above differential equation.
The definition of $G_0^a(u)\approx u^A$ in Theorem \ref{Egypt}
especially implies that each $ J_p^a(W)(z)\in A^a_{\mathrm{Col}}$
for $W\in M'$ is analytic at $z=0$ and $J_p^a(W)(0)=0$
because
$f'\left(G_0^a(u)\cdot u^{-A}\right)=f'\left(G_0^a(u)\right)=
1+\sum\limits_{W\in M':\text{words}} J_p^a(W)(u) \ W $.
Therefore by using 
$J_p^a(W)(0)=(-1)^{dp(W)}Li_W^a(0)=0$,
we obtain inductively the equality
$J_p^a(W)(u)=(-1)^{dp(W)}Li_W^a(u)$.
\qed
\end{pf}

By combining Lemma \ref{Sweden} with Lemma \ref{Netherland},
we get Theorem \ref{Belgium}.(a).

\begin{nota}\label{Denmark}
Let $A^{\land}_{\bold C_p}[[\alpha]]:=A^{\land}_{\bold C_p}\widehat{\otimes}\bold C_p[[\alpha]]$
be the one variable formal power series ring with 
coefficients in
the non-commutative algebra
$A^{\land}_{\bold C_p}$.
Let $g'_1:A^{\land}_{\bold C_p}\to A^{\land}_{\bold C_p}[[\alpha]]$
be the algebra homorphism which sends $A$, $B$ to $A-\alpha$, $B$
respectively and let
$g'_2: A^{\land}_{\bold C_p}[[\alpha]]\to A^{\land}_{\bold C_p}$
be the well-defined $\bold C_p$-linear map which sends
$W\otimes\alpha^q$ to $WA^q$ for each word $W$ and $q\geqslant 0$.
\end{nota}

Consider the $\bold C_p$-linear map
$g'_2\circ g'_1:A^{\land}_{\bold C_p}\to A^{\land}_{\bold C_p} $.

\begin{lem}\label{Austria}
$g'_2\circ g'_1\circ f'=g'_2\circ g'_1$.
\end{lem}

\begin{pf}
By definition, we get easily
$g'_2\circ g'_1(VA)=0$ for $V\in A^{\land}_{\bold C_p} $.
\qed
\end{pf}

\begin{lem}\label{Italy}
$G^a_0(z)=g'_2\circ g'_1\Bigl(f'\bigl(G^a_0(z)\bigr)\Bigr)\cdot z^A$.
\end{lem}

\begin{pf}
By Lemma \ref{Austria}, we get 
\begin{equation}\label{Hotaka}
g'_2\circ g'_1\Bigl(f'\bigl(G^a_0(z)\bigr)\Bigr)=
g'_2\circ g'_1\bigl(G^a_0(z)\bigr)\quad .
\end{equation}
Both $G^a_0(A-\alpha, B)(u)$ and $u^{-\alpha}G^a_0(A,B)(u)$ 
are solutions of the $p$-adic differential equation
$\frac{dG}{du}(u)=(\frac{A-\alpha}{u}+\frac{B}{u-1})G(u)$ in 
$A^a_{\mathrm{Col}}\widehat\otimes A^{\land}_{\bold C_p}[[\alpha]]$
and satisfies the same asymptotic behavior $G(u)\approx u^{A-\alpha}$
as $u\to 0$.
Therefore the uniqueness of solution of the above 
$p$-adic differential equation 
(which can be shown in a way similar to Proposition \ref{Israel}),
we get $G^a_0(A-\alpha, B)(u)=G^a_0(A,B)(u)\cdot u^{-\alpha}$,
from which it follows that 
\begin{equation}\label{Aso}
g'_2\circ g'_1\bigl(G^a_0(z)\bigr)=G^a_0(z)\cdot z^{-A} \quad .
\end{equation}
By \eqref{Hotaka} and \eqref{Aso}, we get 
$G^a_0(z)=g'_2\circ g'_1\Bigl(f'\bigl(G^a_0(z)\bigr)\Bigr)\cdot z^A$.
\qed
\end{pf}

Therefore we see that 
$P^a(z)=g'_2\circ g'_1\Bigl(f'\bigl(G^a_0(z)\bigr)\Bigr)$
(for $P^a(z)$, see the proof of Theorem \ref{Egypt}).

\begin{nota}\label{Lille}
\begin{enumerate}
\renewcommand{\labelenumi}{(\roman{enumi})}
\item
We define the $\bold Q$-bilinear {\sf inner product}
$<\cdot,\cdot>:\bold A\centerdot\times\bold A\centerdot\to\bold A\centerdot$
by $<W,W'>:=\delta_{W,W'}$ for each word (or $1$) $W$ and $W'$,
where 
\[
\delta_{W,W'}:=
\begin{cases}
1, \qquad \text{if } W=W', \\
0, \qquad \text{if } W\neq W', \\
\end{cases}
\]
\item
We define $F':\bold A\centerdot\to\bold A\centerdot$ to be 
the graded $\bold Q$-linear map
which sends each word (or $1$) $W=W'A^r$ ($r\geqslant 0$, $W'\in M'$ or $W'=1$) to
$(-1)^rf'(W'\circ A^r)$.
We note that $Im F'\subseteq M'$.
\end{enumerate}
\end{nota}

\begin{lem}\label{Rennes}
The linear map $F'$ is the transpose of $g'_2\circ g'_1$, i.e.
\[
<(g'_2\circ g'_1)(W_1),W_2>=<W_1,F'(W_2)> \ 
\text{for any words $W$ and $W'$.}
\]
\end{lem}

\begin{pf}
In the case when $W_1\not\in M'$, it is clear.
So we may assume that $W_1\in M'$.
Denote $W_2=W'_2A^r$ ($r\geqslant 0$, $W_2'=1$ or $W'_2\in M'$).
Then, by a direct computation, we deduce elementarily the following
\[
<(g'_2\circ g'_1)(W_1),W_2>=<(g'_2\circ g'_1)(W_1),W'_2A^r>
=<W_1,(-1)^rf'(W_2'\circ A^r)>.
\]
\qed
\end{pf}

{\bf Proof of Theorem \ref{Belgium}.}
By Lemma \ref{Netherland} and Lemma \ref{Rennes}, we get 
\[
g'_2\circ g'_1\Bigl(f'\bigl(G^a_0(z)\bigr)\Bigr)=
1+\sum\limits_{W:\text{words}} J'(W)(z) \ W
\]
where
\begin{align*}
J'(W)(z)&=(-1)^{dp(W)}Li_W^a(z)&
\text{if } W &\in M',\\
J'(W)(z)&=J'(W'A^r)(z)=J'\left(F'(W' A^r)\right)(z)&
= (-1)^r & J'\left(f'(W'\circ A^r)\right)(z)\\
&=(-1)^{dp(W)+r}Li_{f'(W'\circ A^r)}^a(z)&
\text{if } W &=W'A^r \ (r\geqslant 0, W'\in M'),\\
J'(W)(z)&=0&
\text{if } W &=A^r (r\geqslant 1).\\
\end{align*}
By Lemma \ref{Italy}, we obtain
\begin{align*}
G^a_0(z) &=1+\sum\limits_{W:\text{words}} J_p^a(W)(z) \ W \\
& =\left(1+\sum\limits_{W:\text{words}} J'(W)(z) \ W\right)\cdot
\left(\sum\limits_{n=0}^{\infty}\frac{\{log^a(z)\}^n}{n!}A^n\right).
\end{align*}
Then, by a direct calculation, we can show 
the explicit formula in Theorem \ref{Belgium}.
\qed

\begin{nota}\label{Poland}
We denote the involution of $A^{\land}_{\bold C_p}$ which exchanges
$A$ and $B$ by $\tau:A^{\land}_{\bold C_p}\to A^{\land}_{\bold C_p}$.
\end{nota}

By Proposition \ref{England}, we get 

\begin{lem}\label{Hungary}
\[
G^a_1(z)=1+\sum\limits_{W:\text{words}} J_p^a\bigl(\tau(W)\bigr)(1-z)\cdot W\quad .
\]
\end{lem}

\begin{eg}\label{Swiss}
The following is a low degree part of $G^a_0(A,B)(z)$.
\begin{align*}
G^a_0&(A,B)(z)= 1+(log^a z)A+{log^a (1-z)}B  +\frac{(log^a z)^2}{2}A^2-Li^a_2(z)AB\\ 
& +\left\{Li^a_2(z)+(log^a z) log^a(1-z)\right\}BA 
+\frac{\{log^a (1-z)\}^2}{2}B^2
+\frac{(log^a z)^3}{6}A^3\\
&-Li^a_3(z)A^2B
  +\left\{2Li^a_3(z)+(log^a z) Li^a_2(z)\right\}ABA 
+Li^a_{1,2}(z)AB^2 \\
& -\left[Li^a_3(z)-(log^a z) Li^a_2(z)-\frac{(log^a z)^2log^a(1-z)}{2}\right]BA^2 
+Li^a_{2,1}(z)BAB \\
&-\left[Li^a_{1,2}(z)+Li^a_{2,1}(z)-\frac{log^a z\{log^a(1-z)\}^2}{2}\right]B^2A \\
&+\frac{\{log^a (1-z)\}^3}{6}B^3+\cdots \qquad .\\
\end{align*}
\end{eg}

\subsection{Explicit formulae of the $p$-adic Drinfel'd associator}\label{pop}
In this subsection, first we give a proof of Theorem \ref{Chiba}
and then give an explicit formula to express each coefficient of the 
$p$-adic Drinfel'd associator $\Phi^p_{KZ}(A,B)$.
The technique employed here is essentially a $p$-adic analogue of 
that employed in \cite{LM} Appendix A.

\begin{lem}\label{America}
\[
\underset{\epsilon\in\bold C_p}{\underset{\epsilon\to 0}{\lim}^\prime} \epsilon^{-A}
G^a_0(\epsilon)=1 \ .
\]
\end{lem}

\begin{pf}
Since  $P^a(u)=G^a_0(u)\cdot u^{-A}$ lies in 
$A(]0[)\widehat{\otimes}A^{\land}_{\bold C_p}$
and takes value $1$ at $u=0$ by Theorem \ref{Egypt}, we get an expression
$P^a(u)=1+uk(u)$ where $k(u)\in A(]0[)\widehat{\otimes}A^{\land}_{\bold C_p}$.
Thus 
\begin{align*}
\underset{\epsilon\in\bold C_p}{\underset{\epsilon\to 0}{\lim}^\prime} \epsilon^{-A}G^a_0(\epsilon)
=&
\underset{\epsilon\in\bold C_p}{\underset{\epsilon\to 0}{\lim}^\prime} \epsilon^{-A}P^a(\epsilon)\epsilon^A \\
=&
\underset{\epsilon\in\bold C_p}{\underset{\epsilon\to 0}{\lim}^\prime}1+\epsilon\cdot exp\{-log^a(\epsilon)A\}\cdot k(\epsilon)\cdot exp\{log^a(\epsilon)A\} .
\end{align*}
By taking its word expansion and applying Lemma \ref{Russia}
in each term, we get the lemma.
\qed
\end{pf}

Although 
$\underset{\epsilon\in\bold C_p}{\underset{\epsilon\to 0}{\lim}}G^a_0(\epsilon)\epsilon^{-A}=1$ by definition and
$\underset{\epsilon\in\bold C_p}{\underset{\epsilon\to 0}{\lim}^\prime} 
\epsilon^{-A}G^a_0(\epsilon)=1$ 
by Lemma \ref{America},
$\underset{\epsilon\in\bold C_p}{\underset{\epsilon\to 0}{\lim}} \epsilon^{-A}G^a_0(\epsilon)=1$
does not hold.

\begin{lem}\label{Canada}
\[
\underset{\epsilon\in\bold C_p}{\underset{\epsilon\to 0}{\lim}^\prime} \epsilon^{-B}G^a_0(1-\epsilon)=\Phi^p_{KZ}(A,B) \ .
\]
\end{lem}

\begin{pf}
By Proposition \ref{England} and Lemma \ref{America},
we get 
\[
\underset{\epsilon\in\bold C_p}{\underset{\epsilon\to 0}{\lim}^\prime} \epsilon^{-B}G^a_1(1-\epsilon)
=\underset{\epsilon\in\bold C_p}{\underset{\epsilon\to 0}{\lim}^\prime} \epsilon^{-B}G^a_0(B,A)(\epsilon)=1 \ .
\]
Thus
\[
\underset{\epsilon\in\bold C_p}{\underset{\epsilon\to 0}{\lim}^\prime} \epsilon^{-B}G^a_0(1-\epsilon)=
\underset{\epsilon\in\bold C_p}{\underset{\epsilon\to 0}{\lim}^\prime} \epsilon^{-B}G^a_1(1-\epsilon)\Phi^p_{KZ}(A,B)
=\Phi^p_{KZ}(A,B) \ .
\]
\qed
\end{pf}

It is interesting to compare 
$\underset{\epsilon\in\bold C_p}{\underset{\epsilon\to 0}{\lim}^\prime} \epsilon^{-B}G^a_0(1-\epsilon)=\Phi^p_{KZ}(A,B)$
in Lemma \ref{Canada} with 
$\underset{\epsilon\in\bold C_p}{\underset{\epsilon\to 0}{\lim}^\prime} 
\epsilon^{-A}G^a_0(\epsilon)=1$
in Lemma \ref{America}.

{\bf Proof of Theorem \ref{Chiba}.}
By Lemma \ref{Canada}, we obtain 
\[
\underset{\epsilon\in\bold C_p}{\underset{\epsilon\to 0}{\lim}^\prime}
\left(\sum\limits_{n=0}^{\infty}\frac{\{-log^a(\epsilon)\}^n}{n!}B^n
\right)\cdot
\left(1+\sum\limits_{W:\text{words}} J^a_p(W)(1-\epsilon) \ W\right)
=\Phi^p_{KZ}(A,B).
\]
Therefore especially for a word $W\in A\cdot \bold A\centerdot$,
we see that 
$\underset{\epsilon\in\bold C_p}{\underset{\epsilon\to 0}{\lim}^\prime}
J^a_p(W)(1-\epsilon)$
converges to the coefficient of $W$ on $\Phi^p_{KZ}(A,B)$.
Thus for each word $W=A^{k_m-1}B\cdots A^{k_1-1}B$ ($k_i\geqslant 1$)
where $k_m>1$,
we can say that 
$\underset{\epsilon\in\bold C_p}{\underset{\epsilon\to 0}{\lim}^\prime}
J^a_p(W)(1-\epsilon)
=(-1)^m \underset{\epsilon\in\bold C_p}{\underset{\epsilon\to 0}{\lim}^\prime} 
Li^a_{k_1,\cdots,k_m}(1-\epsilon)$
(cf. Theorem \ref{Belgium}.(a)) converges.
\qed

\begin{nota}\label{Mexico}
\begin{enumerate}
\renewcommand{\labelenumi}{(\arabic{enumi})}
\item
Put 
$M=A\cdot\bold A\centerdot\cdot B=\{A\cdot F\cdot B| F\in\bold A\centerdot\}$,
which is a $\bold Q$-linear subspace of 
$\bold A\centerdot$.
Note that there is a natural surjection from
$\bold A\centerdot$  to
$\bold A\centerdot\Bigm/(B\cdot\bold A\centerdot+\bold A\centerdot\cdot A)$.
By identifying the latter space with 
$\bold Q\cdot 1+M (=\bold Q\cdot 1+A\cdot\bold A\centerdot\cdot B)$,
we obtain the $\bold Q$-linear map
$f:\bold A\centerdot\twoheadrightarrow\bold A\centerdot\Bigm/(B\cdot\bold A\centerdot+\bold A\centerdot\cdot A)
\overset{\sim}{\to}\bold Q\cdot 1+M\hookrightarrow\bold A\centerdot$.
Abusively we denote by $f$ the $\bold C_p$-linear map
$ A^{\land}_{\bold C_p}\to A^{\land}_{\bold C_p}$
induced by
$f:\bold A\centerdot\to\bold A\centerdot$.
\item
For each word $A^{p_1}B^{q_1}A^{p_2}B^{q_2}\dotsm A^{p_k}B^{q_k}$
($p_i,q_i\geqslant 1$) in $M$, we define
\begin{align*}
Z_p(W):
&=\underset{\epsilon\in\bold C_p}{\underset{\epsilon\to 0}{\lim}^\prime} Li^a_W(1-\epsilon) \\
&=\zeta_p(\underbrace{1,\ldots 1}_{q_k-1},
p_k+1,
\underbrace{1,\ldots 1}_{q_{k-1}-1},
p_{k-1}+1,\ldots\ldots ,1,p_1+1) \ .
\end{align*}
By extending linearly, we get the $\bold Q$-linear map
$Z_p:M\to\bold C_p$
which sends each word $W$ in $M$ to $Z_p(W)\in\bold C_p$.
\end{enumerate}
\end{nota}

We already know that $Z_p(W)$ is independent of any choice of branch parameter
$a\in\bold C_p$ by Theorem \ref{Gunma}
and lies in $\bold Q_p$ by Theorem \ref{Ishikawa}.

\begin{rem}\label{Cuba}
By combining Lemma \ref{Canada} with Theorem \ref{Ireland}, we get
another proof of branch independency (Theorem \ref{Gunma}) of the value
$\underset{z\in\bold C_p-\{1\}}{\underset{z\to 1}{\lim}^\prime} Li_{k_1,\cdots,k_m}^a(z)$
for $k_m>1$.
\end{rem}

\begin{thm}[Explicit Formulae]\label{Brazil}
Expand the $p$-adic Drinfel'd associator:
$\Phi^p_{KZ}(A,B)=1+\sum\limits_{W:\text{words}} I_p(W)W $.
Then each coefficient $I_p(W)$ can be expressed as follows.
\begin{enumerate}
\renewcommand{\labelenumi}{(\alph{enumi})}
\item  When $W$ is in $M$, \ \ \ $I_p(W)=(-1)^{dp (W)}Z_p(W)$.
\item When $W$ is written as $B^rVA^s (r,s\geqslant 0,V\in M)$,
\[
I_p(W)=(-1)^{dp (W)}\sum_{0\leqslant a\leqslant r,0\leqslant b\leqslant s}(-1)^{a+b}
Z_p\Bigl(f(B^a\circ B^{r-a}VA^{s-b}\circ A^b)\Bigr).
\]
\item When $W$ is written as $B^rA^s (r,s\geqslant 0)$,
\[
I_p(W)=(-1)^{dp (W)}\sum_{0\leqslant a\leqslant r ,0\leqslant  b\leqslant s}(-1)^{a+b}
Z_p\Bigl(f(B^a\circ B^{r-a}A^{s-b}\circ A^b)\Bigr).
\]
\end{enumerate}
\end{thm}

The proof of this theorem will be given in the end of this subsection.

\begin{rem}\label{Peru}
\begin{enumerate}
\renewcommand{\labelenumi}{(\arabic{enumi})}
\item
These explicit formulae are the $p$-adic version of those given
in \cite{F1} Proposition 3.2.3.
\item
Suppose that $k_i\geqslant 1$ and $k_m=1$.
In the complex case, $(-1)^mI(A^{k_m-1}B\cdots A^{k_1-1}B)$
(for $I(\cdot)$, see \cite{F1} Proposition 3.2.3) is called the
{\sf regularized MZV} corresponding to $A^{k_m-1}B\cdots A^{k_1-1}B$
(just something like a modification of the divergent series
$\zeta(k_1,\cdots,k_{m-1},1)$), see for example \cite{IKZ}.
Therefore we may call $(-1)^mI_p(A^{k_m-1}B\cdots A^{k_1-1}B)$
by the {\sf regularized $p$-adic MZV}
corresponding to $A^{k_m-1}B\cdots A^{k_1-1}B$.
\end{enumerate}
\end{rem}

\begin{lem}\label{Syria}
\[
f(\Phi^p_{KZ}(A,B))=1+\sum\limits_{W\in M:\text{words}} I_p(W)W .
\]
\end{lem}

\begin{pf}
Apply $f$ term by term.
\qed
\end{pf}

\begin{lem}\label{Iran}
\[
f(\Phi^p_{KZ}(A,B))=1+\sum\limits_{W\in M:\text{words}}(-1)^{dp(W)}Z_p(W)\cdot W
\]
\end{lem}

\begin{pf}
It follows from Theorem \ref{Belgium}.(a) and
Lemma \ref{Canada}.
\qed
\end{pf}

By combining Lemma \ref{Syria} and Lemma \ref{Iran},
we get Theorem \ref{Brazil}.(a).

\begin{nota}\label{Iraq}
Let $A^{\land}_{\bold C_p}[[\alpha,\beta]]:=A^{\land}_{\bold C_p}\widehat{\otimes}\bold C_p[[\alpha,\beta]]$
be the two variable formal power series ring with 
coefficients in the non-commutative algebra
$A^{\land}_{\bold C_p}$.
Let $g_1:A^{\land}_{\bold C_p}\to A^{\land}_{\bold C_p}[[\alpha,\beta]]$
be the algebra homomorphism which sends $A$, $B$ to $A-\alpha$, $B-\beta$
respectively and let
$g_2: A^{\land}_{\bold C_p}[[\alpha,\beta]]\to A^{\land}_{\bold C_p}$
be the well-defined $\bold C_p$-linear map which sends
$W\otimes\alpha^p\beta^q$ to $B^qWA^p$ for each word $W$ and $p,q\geqslant 0$.
\end{nota}

Consider the $\bold C_p$-linear map 
$g_2\circ g_1:A^{\land}_{\bold C_p}\to A^{\land}_{\bold C_p}$.

\begin{lem}\label{India}
$g_2\circ g_1\circ f=g_2\circ g_1$.
\end{lem}

\begin{pf}
By definition, we get easily 
$g_2\circ g_1(VA)=0$ and $g_2\circ g_1(BV)=0$ for $V\in A^{\land}_{\bold C_p}$.
\qed
\end{pf}

\begin{lem}\label{Thailand}
$\Phi^p_{KZ}(A,B)=g_2\circ g_1\Bigl(f\bigl(\Phi^p_{KZ}(A,B)\bigr)\Bigl)$  .
\end{lem}

\begin{pf}
By Lemma \ref{India}, we get
\begin{equation}\label{Hiei}
g_2\circ g_1\Bigl(f\bigl(\Phi^p_{KZ}(A,B)\bigr)\Bigl)=
g_2\circ g_1\bigl(\Phi^p_{KZ}(A,B)\bigr) \ .
\end{equation}
Both $G^a_0(A-\alpha,B-\beta)(u)$ and $u^{-\alpha}(1-u)^{-\beta}G_0(A,B)(u)$
are solutions of the $p$-adic differential equation 
$\frac{dG}{du}(u)=(\frac{A-\alpha}{u}+\frac{B-\beta}{u-1})G(u)$ in
$A^a_{\mathrm{Col}}\widehat\otimes A^{\land}_{\bold C_p}[[\alpha,\beta]]$
and satisfy the same asymptotic behavior $G(u)\approx u^{A-\alpha}$
as $u\to 0$.
By the uniqueness of solutions of the above $p$-adic differential equation 
(which can be shown in a way similar to Proposition \ref{Israel}), we get
\[
G^a_0(A-\alpha,B-\beta)(u)=u^{-\alpha}(1-u)^{-\beta}G_0(A,B)(u).
\]
Similarly we get 
\[
G^a_1(A-\alpha,B-\beta)(u)=u^{-\alpha}(1-u)^{-\beta}G_1(A,B)(u).
\]
Therefore
\[
G^a_1(A-\alpha,B-\beta)(u)^{-1}G^a_0(A-\alpha,B-\beta)(u)=G^a_1(A,B)(u)^{-1}G^a_0(A,B)(u),
\]
from which it follows that 
\[
\Phi^p_{KZ}(A-\alpha,B-\beta)=\Phi^p_{KZ}(A,B).
\]
Thus we get,
\begin{equation}\label{Shigi}
g_2\circ g_1\bigl(\Phi^p_{KZ}(A,B)\bigr)=\Phi^p_{KZ}(A,B).
\end{equation}
From \eqref{Hiei} and \eqref{Shigi}, it follows that
$\Phi^p_{KZ}(A,B)=g_2\circ g_1\Bigl(f\bigl(\Phi^p_{KZ}(A,B)\bigr)\Bigl)$.
\qed
\end{pf}

\begin{nota}\label{Lyon}
We denote $F:\bold A\centerdot\to\bold A\centerdot$ to be 
the graded $\bold Q$-linear map
which sends each word (or $1$) 
$W=B^rW'A^s$ ($r,s\geqslant 0$, $W'\in M$ or $W'=1$) to
$\underset{0\leqslant a \leqslant r,0\leqslant b \leqslant s}{\sum} (-1)^{a+b}f'(B^a\circ B^{r-a}W'A^{s-b}\circ A^b)$.
We note that $Im F\subseteq M$.
\end{nota}

\begin{lem}\label{Nice}
The linear map $F$ is the transpose of $g_2\circ g_1$.
\end{lem}

\begin{pf}
By an argument similar to Lemma \ref{Rennes}, we can prove this lemma.
\qed
\end{pf}

{\bf Proof of Theorem \ref{Brazil}.}
By combining Lemma \ref{Iran} with Lemma \ref{Thailand} and Lemma \ref{Nice},
we can show the explicit formulae in Theorem \ref{Brazil}
by an argument similar to the proof of Theorem \ref{Belgium}.
\qed

\begin{eg}\label{Laos}
The following is a low degree part of the $p$-adic Drinfel'd associator 
$\Phi^p_{KZ}(A,B)$.
\begin{align*}
\Phi^p_{KZ}& (A,B)=1-\zeta_p(2)AB+\zeta_p(2)BA
-\zeta_p(3)A^2B+2\zeta_p(3)ABA  \\
&+\zeta_p(1,2)AB^2-\zeta_p(3)BA^2-2\zeta_p(1,2)BAB+\zeta_p(1,2)B^2A \\
&-\zeta_p(4)A^3B+3\zeta_p(4)A^2BA+\zeta_p(1,3)A^2B^2-3\zeta_p(4)ABA^2 \\
&+\zeta_p(2,2)ABAB-(2\zeta_p(1,3)+\zeta_p(2,2))AB^2A-\zeta_p(1,1,2)AB^3 \\
&+\zeta_p(4)BA^3-(2\zeta_p(1,3)+\zeta_p(2,2))BA^2B
+(4\zeta_p(1,3)+\zeta_p(2,2))BABA \\
&+3\zeta_p(1,1,2)BAB^2-\zeta_p(1,3)B^2A^2-3\zeta_p(1,1,2)B^2AB \\
&+\zeta_p(1,1,2)B^3A+\cdots .
\end{align*}
\end{eg}

\subsection{Proofs of main results and functional equations among $p$-adic 
multiple polylogarithms}\label{rock}
Here we show functional equations among $p$-adic MPL's in
Theorem \ref{Malaysia} and give proofs of Theorem \ref{Niigata} and
Theorem \ref{Nagano}.

\begin{thm}[Functional Equation among $p$-adic MPL's]\label{Malaysia}
Let $W$ be a word and $z\in\bold P^1(\bold C_p)\backslash\{0,1,\infty\}$, then
\[
J_p^a(W)(1-z)=\underset{W=W'W''}{\sum_{W',W'':\text{words}}}
J_p^a\bigl(\tau(W')\bigr)(z)\cdot I_p(W'') \ .
\]
\end{thm}

For $J_p^a$, $I_p$, $\tau$, see Theorem \ref{Belgium}, Theorem \ref{Brazil}
and Notation \ref{Poland} respectively.
This formulae may be regarded as a functional equation among $p$-adic MPL's
because $J_p^a(\tau(W'))(z)$ (resp. $I_p(W'')$) is expressed in terms of 
$p$-adic MPL's (resp. $p$-adic MZV's).

\begin{pf}
It follows from $G^a_0(A,B)(z)=G^a_1(A,B)(z)\cdot\Phi^p_{KZ}(A,B)$
and Lemma \ref{Hungary}.
\qed
\end{pf}

\begin{eg}\label{Indonesia}
Put $a\in\bold C_p$.
\begin{enumerate}
\renewcommand{\labelenumi}{(\alph{enumi})}
\item
Take $W=BAB$. Then we get 
\[
Li^a_{2,1}(1-z)=2Li^a_3(z)-log^a(z) Li^a_2(z)-\zeta_p(2)log^a(z)
-2\zeta_p(3).
\]
\item 
Take $W=BA^2B$. Then we get 
\begin{align*}
Li^a_{3,1}(1-z)= 
& -2Li^a_{1,3}(z)-Li^a_{2,2}(z)+log^a(z) Li^a_{1,2}(z) \\
&+\zeta_p(2)Li^a_2(z)-\zeta_p(3)log^a(z)-2\zeta_p(1,3)-\zeta_p(2,2).
\end{align*}
\item
Take $W=AB$. Then we get 
\[
Li^a_2(1-z)=-Li^a_2(z)-log^a(z)log^a(1-z)+\zeta_p(2) .
\]
This formula is equal Coleman-Sinnott's functional equation 
of the $p$-adic dilogarithm (\cite{C} Proposition 6.4.(iii))
because $\zeta_p(2)=0$ by Example \ref{Tokyo}.(a).
\end{enumerate}
\end{eg}

{\bf Proof of Theorem \ref{Niigata}.}
By the explicit formulae in Theorem \ref{Belgium} and the functional 
equation of $p$-adic MPL's in Theorem \ref{Malaysia} combined with
Lemma \ref{Russia}, it is immediate to see that
\[
\underset{z\in\bold C_p-\{1\}}{\underset{z\to 1}{\lim}^\prime}J_p^a(W)(1-z)=I_p(W)
\qquad \text{ if it converges.}
\]
Theorem \ref{Niigata} is a special case for 
$W=BA^{k_{m-1}-1}B\cdots A^{k_1-1}B$.
\qed

\begin{nota}\label{Philippines}
\begin{enumerate}
\renewcommand{\labelenumi}{(\arabic{enumi})}
\item
Denote the subset of $A^{\land}_{\bold C_p}$
consisting of formal Lie series in $A^{\land}_{\bold C_p}$
by $\Bbb L^{\land}_{\bold C_p}$
and its topological commutator by 
$[\Bbb L^{\land}_{\bold C_p},\Bbb L^{\land}_{\bold C_p}]$.
We denote $exp \  \Bbb L^{\land}_{\bold C_p}$ 
(resp. $exp \ [\Bbb L^{\land}_{\bold C_p},\Bbb L^{\land}_{\bold C_p}]$ )
to be the subset of $A^{\land}_{\bold C_p}$ consisting of formal power
series in $A^{\land}_{\bold C_p}$
which is an exponential of an element in $\Bbb L^{\land}_{\bold C_p}$
(resp. $[\Bbb L^{\land}_{\bold C_p},\Bbb L^{\land}_{\bold C_p}]$ ).
\item
We define the (non-commutative) $A^a_{\mathrm{Col}}$-algebra homomorphism
\[
\Delta:A^a_{\mathrm{Col}}\underset{\bold C_p}{\widehat\otimes} A^{\land}_{\bold C_p}\to
A^a_{\mathrm{Col}}\underset{\bold C_p}{\widehat\otimes} 
(A^{\land}_{\bold C_p}\underset{\bold C_p}{\widehat\otimes}A^{\land}_{\bold C_p})
\]
to be the homomorphism
which is deduced from $\Delta(A)=A\otimes 1+1\otimes A$ and
$\Delta(B)=B\otimes 1+1\otimes B$.
\end{enumerate}
\end{nota}

\begin{prop}\label{China}
$\Delta(\Phi^p_{KZ})=\Phi^p_{KZ}\widehat\otimes\Phi^p_{KZ}$ .
\end{prop}

\begin{pf}
Put $a\in\bold C_p$.
By the following calculations,
\begin{align*}
&\Delta\bigl(G^a_0(A,B)(u)\bigr)=G^a_0\bigl(\Delta(A),\Delta(B)\bigr)(u)
\approx u^{\Delta(A)} \qquad \text{ as } \ u\to 0, \\
&\frac{d\Delta\bigl(G^a_0(A,B)\bigr)}{du}(u) 
=\frac{dG^a_0\bigl(\Delta(A),\Delta(B)\bigr)}{du}(u) \\
&\qquad\qquad\qquad\qquad =(\frac{\Delta(A)}{u}+\frac{\Delta(B)}{u-1})
G^a_0\bigl(\Delta(A),\Delta(B)\bigr)(u) \\
&\qquad\qquad\qquad\qquad =(\frac{\Delta(A)}{u}+\frac{\Delta(B)}{u-1})\Delta\bigl(G^a_0(A,B)(u)\bigr),\\
\end{align*}
\begin{align*}
&G^a_0(A,B)(u)\widehat\otimes G^a_0(A,B)(u)
=\bigl(G^a_0(A,B)(u)\widehat\otimes 1\bigr)\cdot
\bigl(1\widehat\otimes G^a_0(A,B)(u)1\bigr)\\
&\qquad\qquad\qquad\qquad\qquad\quad
 \approx u^A\widehat\otimes u^A \qquad \text{ as } \ u\to 0, \\
\end{align*}
\begin{align*}
&\frac{d\bigl(G^a_0(A,B)(u)\widehat\otimes G^a_0(A,B)(u)\bigr)}{du}=
\frac{d}{du}\Bigl\{\bigl(G^a_0(A,B)(u)\widehat\otimes 1\bigr)\cdot
\bigl(1\widehat\otimes G^a_0(A,B)(u)\bigr)\Bigr\}    \\
&=\Bigl\{\frac{d}{du}\bigl(G^a_0(A,B)(u)\bigr)\widehat\otimes 1\Bigr\}\cdot
\Bigl\{1\widehat\otimes G^a_0(A,B)(u)\Bigr\}   \\
&\qquad\qquad\qquad
+\Bigl\{ G^a_0(A,B)(u)\widehat\otimes 1\Bigr\}\cdot
\Bigl\{1\widehat\otimes\frac{d}{du}\bigl(G^a_0(A,B)(u)\bigr)\Bigr\}  \\
&=\Bigl\{\bigl(\frac{A}{u}+\frac{B}{u-1})\cdot G^a_0(A,B)(u)\widehat\otimes 1\Bigr\}\cdot
\Bigl\{1\widehat\otimes G^a_0(A,B)(u)\Bigr\}   \\
&\qquad\qquad\qquad
+\Bigl\{ G^a_0(A,B)(u)\widehat\otimes 1\Bigr\}\cdot
\Bigl\{1\widehat\otimes\bigl(\frac{A}{u}+\frac{B}{u-1})\cdot G^a_0(A,B)(u)\Bigr\}   \\
&=\bigl(\frac{A\widehat\otimes 1 +1 \widehat\otimes A}{u}+\frac{B\widehat\otimes 1 +1 \widehat\otimes 1}{u-1})\cdot
\Bigl\{ G^a_0(A,B)(u)\widehat\otimes G^a_0(A,B)(u)\Bigr\},
\end{align*}
we see that both $\Delta\bigl(G^a_0(A,B)(u)\bigr)$ and
$G^a_0(A,B)(u)\widehat\otimes G^a_0(A,B)(u)$
are solutions in 
$A^a_{\mathrm{Col}}\underset{\bold C_p}{\widehat\otimes}(A^{\land}_{\bold C_p}\underset{\bold C_p}{\widehat\otimes}A^{\land}_{\bold C_p})$
of the $p$-adic differential equation
\[
\frac{dH}{dt}(t)=(\frac{\Delta(A)}{t}+\frac{\Delta(B)}{t-1})\cdot H(t)
\]
which satisfies $H(t)\approx t^A \widehat\otimes t^A$ as $t\to 0$.
Because of the uniqueness of solution for above $p$-adic differential equation
(which can be shown in a similar way to Proposition \ref{Israel}), we get
\[
\Delta\bigl(G^a_0(A,B)(u)\bigr)=G^a_0(A,B)(u)\widehat\otimes G^a_0(A,B)(u).
\]
Similarly we get 
\[
\Delta\bigl(G^a_1(A,B)(u)\bigr)=G^a_1(A,B)(u)\widehat\otimes G^a_1(A,B)(u).
\]
Therefore
\begin{align*}
\Delta\bigl(\Phi^p_{KZ}\bigr)&
=\Delta\Bigl(G^a_1(A,B)(u)^{-1}\cdot G^a_0(A,B)(u)\Bigr)  \\
&=\Delta\Bigl(G^a_1(A,B)(u)\Bigr)^{-1} \cdot
\Delta\Bigl(G^a_0(A,B)(u)\Bigr) \\
&=\Bigl(G^a_1(A,B)(u)\widehat\otimes G^a_1(A,B)(u)\Bigr)^{-1} \cdot
\Bigl(G^a_0(A,B)(u)\widehat\otimes G^a_0(A,B)(u)\Bigr)   \\
&=\Bigl(G^a_1(A,B)(u)^{-1}\cdot G^a_0(A,B)(u)\Bigr)\widehat\otimes
\Bigl(G^a_1(A,B)(u)^{-1}\cdot G^a_0(A,B)(u)\Bigr)  \\
&=\Phi^p_{KZ}\widehat\otimes\Phi^p_{KZ}\qquad .
\end{align*}
\qed
\end{pf} 

\begin{lem}\label{Korea}
Put $a\in\bold C_p$.
Denote $g_0^a(\alpha,\beta)(z)$ and $g_1^a(\alpha,\beta)(z)$
to be the images of $G^a_0(A,B)(z)$ and $G^a_1(A,B)(z)$ 
by the natural projection 
$A^a_{\mathrm{Col}}[[A,B]]\twoheadrightarrow A^a_{\mathrm{Col}}[[\alpha,\beta]]$
sending $A$ (resp. $B$) to $\alpha$ (resp. $\beta$),
where $A^a_{\mathrm{Col}}[[\alpha,\beta]]$
is the two variable commutative formal power series ring with
$A^a_{\mathrm{Col}}$-coefficients.
Then 
\[
g_0^a(\alpha,\beta)(z)=g_1^a(\alpha,\beta)(z)=z^\alpha(1-z)^\beta
\quad\text{ in }\quad
A^a_{\mathrm{Col}}[[\alpha,\beta]] \ .
\]
\end{lem}

\begin{pf}
Both $g_0^a(\alpha,\beta)(z)$ and $z^\alpha(1-z)^\beta$ are solutions in 
$A^a_{\mathrm{Col}}[[\alpha,\beta]]$
of the $p$-adic differential equation 
\[
\frac{dH}{dt}(t)=(\frac{\alpha}{t}+\frac{\beta}{t-1})\cdot H(t)
\]
which satisfy $H(t)\approx t^\alpha$ as $t\to 0$.
By the uniqueness of solution of the above $p$-adic differential equation
(which can be shown in a similar way to Proposition \ref{Israel}), we get
\[
g_0^a(\alpha,\beta)(z)=z^\alpha(1-z)^\beta .
\]
Similarly we get
\[
g_1^a(\alpha,\beta)(z)=z^\alpha(1-z)^\beta . 
\]
Therefore
\[
g_0^a(\alpha,\beta)(z)=g_1^a(\alpha,\beta)(z)=z^\alpha(1-z)^\beta .
\qquad\qquad\qquad\qquad  \qed
\]
\end{pf}

\begin{thm}\label{Japan}
$\Phi^p_{KZ}(A,B)\in exp \ [\Bbb L^{\land}_{\bold C_p},\Bbb L^{\land}_{\bold C_p}]$.
\end{thm}

\begin{pf}
By Proposition \ref{China}, we see that $\Phi^p_{KZ}$  is group-like,
which means that $\Phi^p_{KZ}(A,B)\in exp \ \Bbb L^{\land}_{\bold C_p}$
(see \cite{Se} Part I Ch IV \S7).
It follows from Lemma \ref{Korea} that
\[
\Phi^p_{KZ}(A,B)= G^a_1(A,B)(z)^{-1}G^a_0(A,B)(z)\in exp \ [\Bbb L^{\land}_{\bold C_p},\Bbb L^{\land}_{\bold C_p}] \ . \quad
\qed
\]
\end{pf}

\begin{cor}[Shuffle Product Formulae]\label{Singapore}  
For each word $W$ and $W'\in M$,
\[
Z_p(W)\cdot Z_p(W')=Z_p(W\circ W') \ .
\]
\end{cor}

\begin{pf}
Consider the graded $\bold Q$-linear map 
$I_p:\bold A\centerdot\to Z^{(p)}_\centerdot$
which sends each word $W$ to $I_p(W)\in Z^{(p)}_\centerdot$
(for $I_p(W)$, see Theorem \ref{Brazil}).
Then by Proposition \ref{China}, we obtain the shuffle product formulae
$I_p(W)\cdot I_p(W')=I_p(W\circ W')$ on $\bold Q_p$
for each word $W'$ and $W''$.
By Theorem \ref{Brazil}, we get the corollary.
\qed
\end{pf}

{\bf Proof of Theorem \ref{Nagano}.}
It follows immediately from Corollary \ref{Singapore}.
\qed

\ifx\undefined\bysame
\newcommand{\bysame}{\leavevmode\hbox to3em{\hrulefill}\,}
\fi


\begin{thebibliography}{Utah}
\bibitem[AK]{AK}
Arakawa, T. and Kaneko, M.;
On poly-Bernoulli numbers,
Comment. Math. Univ. St. Paul. 48 (1999), no. 2, 159--167.

\bibitem[Ber]{Ber}
Berthelot, P.;
Cohomologie rigide et cohomologie rigide \`{a} support propre, 
Premi\`{e}re partie,                            
Pr\'{e}publication IRMAR 96-03, 89 pages (1996).

\bibitem[Bes1]{Be1}
Besser, A.;
Syntomic regulators and p-adic integration II: K$_2$ of curves, 
Israel Journal of Math. 120 (2000), 335-360. 

\bibitem[Bes2]{Be2}
\bysame
Coleman integration using the Tannakian formalism, 
Mathematische Annalen 322 (2002) 1, 19-48. 

\bibitem[BdJ]{BdJ}
Besser, A. and de Jeu, R.;
The syntomic regulator for K-theory of fields,
to be appeared in Annales Scientifiques de L'ecole normale superieure,
available at {\tt arXiv:math.AG/0005069}.

\bibitem[BF]{BF}
Besser, A. and Furusho, H.;
The double shuffle relation for $p$-adic multiple zeta values,
in preparation.

\bibitem[Br]{Br}
Breuil, C.;
Int\'{e}gration sur les vari\'{e}t\'{e}s $p$-adiques 
(d'apr\'{e}s Coleman, Colmez),
S\'{e}minaire Bourbaki, Vol. 1998/99, 
Ast\'{e}risque No. 266, (2000), Exp. No. 860, 5, 319--350.

\bibitem[C]{C}
Coleman, R.;
Dilogarithms, regulators and $p$-adic $L$-functions,
Invent. Math. 69 (1982), no. 2, 171--208. 

\bibitem[CdS]{CdS}
Coleman, R. and de Shalit, E.;
$p$-adic regulators on curves and special values of $p$-adic $L$-functions,
Invent. Math. 93 (1988),no. 2, 239--266. 

\bibitem[Dr]{Dr}
Drinfel'd, V. G.;
On quasitriangular quasi-Hopf algebras and  a group closely connected with ${\rm Gal}(\overline{Q}/{Q})$,
Leningrad Math. J. 2 (1991), no. 4, 829--860. 

\bibitem[F0]{F0}
Furusho, H.;
The multiple zeta value algebra and the stable derivation algebra,
to appear in Publ. Res. Inst. Math. Sci. Vol 39. no 4, 
also available at {\tt arXiv:math.NT/0011261}.

\bibitem[F1]{F1}
\bysame
Multiple zeta values and Grothendieck-Teichm\"{u}ller groups,
RIMS-1357, preprint (2002), submitted.

\bibitem[F2]{F2}
\bysame
$p$-adic multiple zeta values II
--- various realizations of motivic fundamental groups of the projective line
minus three points, in preparation.

\bibitem[F3]{F3}
\bysame
$p$-adic multiple zeta values III
--- Grothendieck-Teichm\"{u}ller groups, in preparation.

\bibitem[Gon]{Gon}
Goncharov, A. B.;
Multiple polylogarithms, cyclotomy and modular complexes,
Math. Res. Lett. 5 (1998), no. 4, 497--516.


\bibitem[H]{H}
Hoffman, M.;
The algebra of multiple harmonic series,
J. Algebra 194 (1997), no. 2, 477--495. 

\bibitem[IKZ]{IKZ}
Ihara, K., Kaneko, M. and Zagier, D.;
Derivation and double shuffle relations for multiple zeta values,
preprint (2002).

\bibitem[Kan]{Kan}
Kaneko, M.;
Poly-Bernoulli numbers,
J. Th\'{e}or. Nombres Bordeaux 9 (1997), no. 1, 221--228. 

\bibitem[Kas]{Kas}
Kassel, C.;
Quantum groups, 
Graduate Texts in Mathematics, 155. Springer-Verlag, New York, 1995. 

\bibitem[KNQ]{KNQ}
Kolster, M. and  Nguyen Quang Do, T.;
Syntomic regulators and special values of $p$-adic $L$-functions,
Invent. Math. 133 (1998), no. 2, 417--447. 

\bibitem[LM]{LM}
Le, T.T.Q. and Murakami, J.;
Kontsevich's integral for the Kauffman polynomial,
Nagoya Math. J. 142 (1996), 39--65.

\bibitem[Se]{Se}
Serre, J-P.;
Lie algebras and Lie groups,
lectures given at Harvard University (1964), Second edition,
Lecture Notes in Mathematics, 1500, Springer-Verlag, Berlin, 1992.

\bibitem[Sou]{Sou}
Soul\'{e}, C.;
On higher $p$-adic regulators, 
Algebraic $K$-theory, Evanston 1980 
(Proc. Conf., Northwestern Univ., Evanston, Ill., 1980),
pp. 372--401, Lecture Notes in Math., 854, Springer, Berlin-New York, 1981.

\end{thebibliography}
\end{document}